\documentclass[a4paper,11pt,reqno]{amsart}
\usepackage[T1]{fontenc}
\usepackage{lmodern}
\usepackage{amsmath,amssymb,mathtools}
\usepackage{microtype}
\usepackage[margin=28mm]{geometry}
\usepackage{enumitem}
\usepackage{amsaddr}
\usepackage{hyperref}
\hypersetup{colorlinks=true, linkcolor=blue, citecolor=blue}
\newtheorem{theorem}{Theorem}[section]
\newtheorem{proposition}[theorem]{Proposition}
\newtheorem{lemma}[theorem]{Lemma}
\newtheorem{corollary}[theorem]{Corollary}
\newtheorem{assumption}[theorem]{Assumption}
\theoremstyle{remark}
\newtheorem{remark}[theorem]{Remark}
\numberwithin{equation}{section}

\newcommand{\R}{\mathbb R}
\newcommand{\eps}{\varepsilon}
\newcommand{\calR}{\mathcal R}
\newcommand{\Tr}{\operatorname{Tr}}
\newcommand{\Lip}{\operatorname{Lip}}
\newcommand{\dist}{\operatorname{dist}}

\newcommand{\one}{\mathbf 1}
\newcommand{\dd}{\,dx}
\newcommand{\norm}[1]{\left\lVert#1\right\rVert}
\newcommand{\abs}[1]{\left\lvert#1\right\rvert}
\setlist[enumerate]{label=\textup{(\roman*)},leftmargin=2.2em,itemsep=3pt}
\title[Smooth exact-segregation approximation]{Gamma convergence and smooth exact-segregation approximation for partially segregated systems}
\author{Farid Bozorgnia}
\address{\vspace{-6pt} Department of Mathematics, New Uzbekistan University,\\ Tashkent, Uzbekistan\\
\texttt{f.bozorgnia@newuu.uz}}

\author{Avetik Arakelyan}
\address{\vspace{-6pt} Institute of Mathematics SNCO, Yerevan, Armenia,\\ Yerevan State University, Yerevan, Armenia,\\ DiGREA Lab, CSIE, Yerevan, Armenia\\
\texttt{avetik.arakelyan@ysu.am}}

\date{}

\begin{document}

\begin{abstract}
\vspace{15pt}
We study three-component Dirichlet energies subject to the partial segregation constraint \(u_1u_2u_3=0\). Our main result is a smooth approximation theorem on bounded Lipschitz domains in arbitrary dimension. For compatible nonnegative Lipschitz boundary data, every admissible \(H^1\) configuration can be approximated strongly in \(H^1\) by globally Lipschitz maps that are smooth in the interior and preserve both the prescribed boundary trace and the segregation constraint exactly. The construction combines a Lipschitz retraction, variable-radius mollification, and smooth truncation preserving existing zeros, without assumptions on interfaces or junctions. Under the weaker assumption of compatible nonnegative \(H^{1/2}\) boundary data, we establish existence of minimizers and \(\Gamma\)-convergence in \(L^2\) for the energies with product penalization \(\varepsilon^{-1}\int_\Omega(u_1u_2u_3)^2\,dx\). We also identify the exact \(L^2\) distance to the constrained class and construct segregated corrections with controlled Dirichlet energy. The resulting correction estimate has order \(\varepsilon^{1/6}\) for bounded-energy configurations; this exponent is sharp in general and improves to \(o(\varepsilon^{1/6})\) for almost-minimizers.
\end{abstract}

\maketitle

\section{Introduction}

Segregation phenomena arise in physical, biological, and ecological systems
in which competing components share a spatial domain while interactions
restrict their coexistence. In elliptic and reaction--diffusion models,
strong competition often leads to configurations governed by Dirichlet-type
energies subject to exclusion constraints. Understanding these limits
requires both a variational description of the competition process and
approximation procedures that respect the resulting constraints.

In the classical \emph{pairwise segregation} problem, the positivity sets
of distinct components are disjoint. This regime appears in strongly
competing Lotka--Volterra systems and in models of competing densities.
Contributions include those of Conti, Terracini, and Verzini
\cite{ContiTerraciniVerzini2005}, Caffarelli and Lin
\cite{CaffarelliLin2007}, Chang et al.\ \cite{chang2004segregated},
and Tavares and Terracini \cite{tavares2003segregated}; see also
\cite{ArakelyanBozorgnia2017}. In \emph{long-range segregation},
competition acts at a distance and imposes a different geometry on the
limiting phases; see
\cite{CaffarelliPatriziQuitalo2017,Bozorgnia2016Acta}.

We consider \emph{partial}, or \emph{ternary}, segregation of three
nonnegative components. For a configuration
\[
U=(u_1,u_2,u_3),
\]
the exclusion constraint is
\[
u_1u_2u_3=0\quad\text{a.e. in }\Omega.
\]
Thus at least one component vanishes at almost every point, while two
components may coexist. This permits overlapping positivity sets and
introduces interfaces between different regions of coexistence.
Singularly perturbed systems leading to this constraint are studied in
\cite{bozorgnia2018singularly}.

The regularity theory for partial segregation has been developed by
Soave and Terracini \cite{soave2024partial1,soave2024partial2}.
Under uniform boundedness and local minimality assumptions, their results
give uniform local $C^{0,\alpha}$ estimates for penalized configurations
for every $0<\alpha<3/4$, and optimal local $C^{0,3/4}$ regularity of the
limits, together with a description of the free boundary.
Interface distributions for constrained minimizers are investigated in
\cite{Bozorgnia2026Interfaces}.

The approximation problem addressed here concerns every admissible
Sobolev configuration, without requiring minimality or regularity of its
interfaces. Standard Sobolev approximation does not by itself preserve
the nonlinear exclusion constraint: componentwise mollification can
produce a region where all three components are positive. A Lipschitz
retraction restores segregation, but its composition with a smooth map
need not remain smooth. Preserving the prescribed boundary trace creates
an additional requirement. Our main result provides approximations that
satisfy all these conditions simultaneously, in arbitrary dimension.

Let $\Omega\subset\R^d$ be a bounded Lipschitz domain, $d\ge1$.
For the variational results, we assume compatible nonnegative boundary
data
\[
\varphi=(\varphi_1,\varphi_2,\varphi_3)
\in H^{1/2}(\partial\Omega)^3,
\qquad
\varphi_i\ge0,\quad \varphi_1\varphi_2\varphi_3=0
\quad\text{a.e. on }\partial\Omega.
\]
Define
\begin{align}
\mathcal A
&=\left\{U=(u_1,u_2,u_3)\in H^1(\Omega)^3:
u_i\ge0\ \text{a.e.},\quad \Tr U=\varphi\right\},
\label{eq:A}\\
S
&=\left\{U\in\mathcal A:
u_1u_2u_3=0\ \text{a.e. in }\Omega\right\},
\label{eq:S}
\end{align}
and
\begin{equation}\label{eq:dirichlet}
E(U)=\int_\Omega\sum_{i=1}^3|\nabla u_i|^2\dd.
\end{equation}
The constrained and penalized energies are defined on $L^2(\Omega)^3$ by
\begin{equation}\label{eq:Ezero}
E^0(U)=
\begin{cases}
E(U),&U\in S,\\
+\infty,&U\notin S,
\end{cases}
\end{equation}
and, for $\eps>0$,
\begin{equation}\label{eq:Eepsilon}
E^\eps(U)=
\begin{cases}
\displaystyle E(U)+\frac1\eps\int_\Omega(u_1u_2u_3)^2\dd,
&U\in\mathcal A,\\[2mm]
+\infty,&U\notin\mathcal A.
\end{cases}
\end{equation}
The penalization integral is allowed to be infinite. This is relevant
for $d\ge4$, where membership in $H^1(\Omega)^3$ does not ensure that
$(u_1u_2u_3)^2$ is integrable. The extended-valued formulation therefore
requires no additional integrability hypothesis depending on the
dimension.

Compatibility of the trace implies $S\ne\varnothing$, and both energies
admit minimizers. We prove equicoercivity and
\[
E^\eps\xrightarrow{\Gamma}E^0
\quad\text{in }L^2(\Omega)^3
\quad\text{as }\eps\to0.
\]
Penalized almost-minimizers approach the set of constrained minimizers
strongly in $H^1$, and their penalization energy tends to zero.
The variational argument uses Sobolev compactness and Fatou's lemma for
the liminf inequality, while the constant sequence gives the limsup
inequality; see \cite{braides2002gamma} for the general framework.

 The constant recovery sequence $U^\eps\equiv U$ carries no more regularity
than $U$ itself, so $\Gamma$-convergence alone does not show that an
arbitrary element of $S$ can be approximated within $S$ by regular
configurations. To obtain recovery sequences consisting of maps in
$\Lip(\overline\Omega)^3\cap C^\infty(\Omega)^3$ with exact trace and exact
segregation, we prove strong density within $S$ under Lipschitz boundary
data. The following theorem is the main result of the paper.

\begin{theorem}[Smooth exact-segregation approximation]\label{thm:main}
Let $\Omega\subset\R^d$ be bounded and Lipschitz, $d\ge1$.
Assume the boundary data satisfy
\begin{equation}\label{eq:lip-data}
\varphi\in\Lip(\partial\Omega)^3,
\quad \varphi_i\ge0,
\quad \varphi_1\varphi_2\varphi_3=0
\quad\text{on }\partial\Omega.
\end{equation}
Then, for every $U\in S$, there exists a sequence
\begin{equation}\label{eq:main-regularity}
U^n\in S\cap\Lip(\overline\Omega)^3\cap C^\infty(\Omega)^3
\end{equation}
such that
\begin{equation}\label{eq:main-convergence}
U^n\to U
\quad\text{strongly in }H^1(\Omega)^3.
\end{equation}
In particular, each approximant satisfies the exact geometric and
boundary constraints:
\[
u_i^n\ge0,\quad u_1^nu_2^nu_3^n\equiv0\ \text{in }\Omega,
\quad U^n|_{\partial\Omega}=\varphi.
\]
Furthermore, for every sequence $\eps_n\to0$,
\begin{equation}\label{eq:main-energy}
E^{\eps_n}(U^n)=E^0(U^n)\to E^0(U).
\end{equation}
\end{theorem}

Here $\Lip$ denotes Lipschitz continuity with respect to Euclidean
distance. Smoothness is asserted in the interior, with the prescribed
Lipschitz values remained on the boundary. The theorem applies to
arbitrary elements of $S$ and requires no decomposition or regularity
assumption on their interfaces. It also implies equality of the energy
infima over the Sobolev and regular admissible classes.
The approximation is qualitative: neither a uniform Lipschitz bound
along the sequence nor a prescribed relation between the smoothing
scale and $\eps_n$ is asserted.
 
The proof consists of two stages. First, an unconstrained approximation
with fixed trace is mapped into $S$ by a globally Lipschitz retraction
$\calR$ (defined in  \eqref{eq:retraction}), obtained by subtracting the smallest component from each
coordinate. Strong continuity of this retraction in $H^1$ gives
Lipschitz approximations with exact trace and segregation. Second, each
Lipschitz configuration is mollified at a spatially varying radius that
vanishes at the boundary. A smooth truncation, with threshold determined
by the mollification error, preserves every existing zero of this
intermediate configuration. It therefore maintains segregation while
giving interior smoothness and strong $H^1$ convergence. Additional
zeros may be created; preservation of the positivity sets of the
original Sobolev configuration is not required.

The same retraction gives quantitative corrections of penalized
configurations. Under the compatible $H^{1/2}$ trace assumptions,
Section~\ref{sec:correction} proves that
\[
E^\eps(V)\le M
\quad\to\quad
\norm{V-\calR(V)}_{L^2(\Omega)^3}
\le\sqrt3\,|\Omega|^{1/3}(\eps M)^{1/6},
\qquad \calR(V)\in S,
\]
with corrected energy satisfying $E^0(\calR(V))\le3E(V)$.
Proposition~\ref{prop:exact-distance} also identifies the exact distance
to the constrained class
\[
\dist_{L^2(\Omega)^3}(V,S)
=\|\min_i v_i\|_{L^2(\Omega)}
=\frac1{\sqrt3}\|V-\calR(V)\|_{L^2(\Omega)^3}
\]
for every nonnegative $V\in L^2(\Omega)^3$. This distance need not be
attained by a Sobolev competitor. The retraction supplies an admissible
correction together with a bound on its Dirichlet energy.
A family with homogeneous boundary data shows that the exponent $1/6$
is optimal for general bounded-energy configurations
(Proposition~\ref{prop:sharpness}). For almost-minimizers, the vanishing
penalization energy improves the correction estimate to
$o(\eps^{1/6})$ (Corollary~\ref{cor:little-o}).

Numerical treatments of pairwise segregation include
\cite{ArakelyanBozorgnia2012,Bozorgnia2009,arakelyan2016numerical,arakelyan2018convergence},
and algorithms for the ternary constraint are developed in \cite{bn}.
Related competition models are considered in
\cite{Vauchelet2025,Xu2014,LanzaraMontefusco2019}.

The paper is organized as follows: 
Section~\ref{sec:minproblem} recalls the relevant regularity results.
Section~\ref{sec:retraction} establishes the properties of the retraction,
admissibility, and existence of minimizers.
Section~\ref{sec:gamma} treats variational convergence, and
Section~\ref{sec:density} proves Theorem~\ref{thm:main}.
Section~\ref{sec:correction} develops the quantitative correction and
distance estimates.

\section{The minimization problem and known results}\label{sec:minproblem}

We use the classes $\mathcal A$, $S$ and the energies introduced above.
Throughout the paper, $\Omega\subset\R^d$ is a bounded Lipschitz domain,
and the standing boundary assumptions are
\[
\varphi\in H^{1/2}(\partial\Omega)^3,
\qquad \varphi_i\ge0,
\qquad \varphi_1\varphi_2\varphi_3=0
\quad\text{a.e. on }\partial\Omega.
\]
Additional hypotheses will be stated where they are needed.
In particular, bounded boundary data are used below to obtain the
penalized equations with bounded coefficients, whereas Lipschitz
boundary data are required for the smooth approximation theorem.

\subsection{Variational formulation and energy bounds}

The constrained and penalized minimization problems are
\[
\min_{U\in S}E(U),
\qquad
\min_{U\in\mathcal A}E^\eps(U),
\quad \eps>0.
\]
Proposition~\ref{prop:nonempty} proves that compatibility of the trace
implies $S\ne\varnothing$, and Proposition~\ref{prop:existence}
establishes existence of minimizers for both problems. Thus nonemptiness
is a consequence of the boundary assumptions. Since the segregation
constraint is nonconvex, strict convexity of the Dirichlet energy on
the affine space of maps with fixed trace does not imply uniqueness of
the constrained minimizer.

Every $W\in S$ is a common finite-energy competitor for all penalty
parameters, since
\[
E^\eps(W)=E^0(W)=E(W).
\]
Consequently, if $U^\eps$ is a global minimizer of $E^\eps$, then
\begin{equation}\label{eq:background-energy-bound}
E(U^\eps)
+\frac1\eps\int_\Omega
\bigl(u_1^\eps u_2^\eps u_3^\eps\bigr)^2\dd
\le E(W).
\end{equation}
Because $U^\eps-W\in H_0^1(\Omega)^3$, Poincar\'e's inequality gives
\[
\sup_{0<\eps<1}\norm{U^\eps}_{H^1(\Omega)^3}<\infty,
\qquad
\norm{u_1^\eps u_2^\eps u_3^\eps}_{L^2(\Omega)}
\le\eps^{1/2}E(W)^{1/2}.
\]
These estimates explain the emergence of segregation: after extracting
an almost everywhere convergent subsequence, Fatou's lemma forces the
product of the limiting components to vanish.
They do not, by themselves, imply that the scaled penalization energy
in \eqref{eq:background-energy-bound} tends to zero. That stronger
conclusion follows from the local theory recalled below and, for global
almost-minimizers under the standing trace assumptions, from the
variational argument in Section~\ref{sec:gamma}.

\subsection{The penalized Euler--Lagrange system}

Suppose in addition that $\varphi_i\in L^\infty(\partial\Omega)$.
Proposition~\ref{prop:bounded-minimizers} shows that every global
minimizer of either problem satisfies
\begin{equation}\label{eq:background-maximum-bound}
0\le u_i\le\norm{\varphi_i}_{L^\infty(\partial\Omega)}
\quad\text{a.e. in }\Omega,
\quad i=1,2,3.
\end{equation}
The bound follows by truncating each component at the corresponding
boundary maximum; this preserves admissibility and cannot increase
either energy.

Although $\mathcal A$ contains only nonnegative maps, the penalized
minimizer satisfies an equation throughout $\Omega$. Indeed, taking
componentwise absolute values preserves the energy integral and the
nonnegative boundary trace. Minimality therefore also holds for the
same integral among signed competitors with trace $\varphi$.
For $\zeta\in C_c^\infty(\Omega)$, variation of the $i$th component
gives
\[
0
=2\int_\Omega\nabla u_i^\eps\cdot\nabla\zeta\dd
+\frac2\eps\int_\Omega
u_i^\eps\prod_{j\ne i}(u_j^\eps)^2\zeta\dd.
\]
The boundedness in \eqref{eq:background-maximum-bound} justifies this
differentiation. Dividing by $2$ yields the weak formulation of
\begin{equation}\label{eq:penalized_system}
\left\{
\begin{aligned}
\Delta u_i^\eps
&=\frac1\eps u_i^\eps\prod_{j\ne i}(u_j^\eps)^2
&&\text{in }\Omega,\\
u_i^\eps&\ge0&&\text{in }\Omega,\\
\Tr u_i^\eps&=\varphi_i&&\text{on }\partial\Omega,
\end{aligned}\right.
\quad i=1,2,3.
\end{equation}
For fixed $\eps>0$, interior elliptic estimates and bootstrapping imply
$U^\eps\in C^\infty(\Omega)^3$.
Each component is subharmonic and, on every connected component of
$\Omega$, is either identically zero or strictly positive.
If one component vanishes identically, the remaining equations reduce
to the Laplace equation. These assertions are proved in
Proposition~\ref{prop:bounded-minimizers}.

\subsection{Known interior regularity and limiting configurations}

The smoothness of $U^\eps$ for fixed $\eps$ does not give derivative
bounds independent of $\eps$, since the reaction term contains the
factor $\eps^{-1}$. The local estimates of Soave and Terracini provide
uniform control through this singular limit. We recall their hypotheses
and the conclusions relevant to the present setting.

\begin{assumption}\label{ass:h1-h2}
Let $d\ge2$, and let
\[
\{u^\eps=(u_1^\eps,u_2^\eps,u_3^\eps)\}_{0<\eps<1}
\subset H^1_{\mathrm{loc}}(\Omega)^3
\]
be a family of positive weak solutions of the interior equations in
\eqref{eq:penalized_system}. Assume:
\begin{itemize}
\item[(h1)] The family is uniformly bounded:
\[
\sup_{0<\eps<1}\norm{u^\eps}_{L^\infty(\Omega)^3}<\infty.
\]
\item[(h2)] Each $u^\eps$ minimizes the energy with respect to
compactly supported variations: for every $\omega\Subset\Omega$ and
every $\Psi\in H_0^1(\omega)^3$,
\[
\begin{aligned}
\int_\omega\left(
\sum_i|\nabla u_i^\eps|^2
+\eps^{-1}\prod_i(u_i^\eps)^2\right)\dd
\le
\int_\omega\left(
\sum_i|\nabla(u_i^\eps+\Psi_i)|^2
+\eps^{-1}\prod_i(u_i^\eps+\Psi_i)^2\right)\dd.
\end{aligned}
\]
\end{itemize}
\end{assumption}

For global penalized minimizers with bounded boundary data,
\eqref{eq:background-maximum-bound} gives (h1), while global
minimality and the signed-competitor argument give (h2).
On connected components where all three boundary traces are nontrivial,
the strong maximum principle gives positivity of all three components,
as required above. If a component of the boundary data, $\varphi_i$, vanishes
identically on the boundary of a connected component of $\Omega$,
the corresponding $u_i^\eps$ vanishes identically in that connected
component, and the other components are harmonic. This case is covered by standard interior
estimates for harmonic functions.

\begin{theorem}[{\cite[Theorem~1.1]{soave2024partial1},
\cite[Theorem~1.6]{soave2024partial2}}]\label{thm:ST}
Under Assumption~\ref{ass:h1-h2}, for every compact set
$K\Subset\Omega$ and every $0<\alpha<3/4$,
\[
\sup_{0<\eps<1}\norm{u^\eps}_{C^{0,\alpha}(K)^3}<\infty.
\]
For every sequence $\eps_n\to0$, there exist a subsequence,
not relabeled, and a nonnegative configuration
$\widetilde U=(\widetilde u_1,\widetilde u_2,\widetilde u_3)$ such that
\[
u^{\eps_n}\to\widetilde U
\quad\text{strongly in }H^1_{\mathrm{loc}}(\Omega)^3
\text{ and }C^{0,\alpha}_{\mathrm{loc}}(\Omega)^3
\quad\text{for every }0<\alpha<3/4.
\]
Moreover,
\[
\widetilde u_1\widetilde u_2\widetilde u_3=0
\quad\text{in }\Omega,
\qquad
\frac1{\eps_n}\int_\omega
\prod_i(u_i^{\eps_n})^2\dd\to0
\quad(\omega\Subset\Omega).
\]
The limit belongs to $C^{0,3/4}_{\mathrm{loc}}(\Omega)^3$.
The exponent $3/4$ is optimal for constrained minimizers, as shown by
the planar example in \cite[Theorem~1.7]{soave2024partial2}.
\end{theorem}

The limiting components also satisfy
\begin{equation}\label{eq:background-limit-harmonicity}
\Delta\widetilde u_i\ge0
\quad\text{in }\mathcal D'(\Omega),
\qquad
\Delta\widetilde u_i=0
\quad\text{in }\{\widetilde u_i>0\},
\quad i=1,2,3.
\end{equation}
Subharmonicity follows by passing to the limit in
\[
\int_\Omega\nabla u_i^{\eps_n}\cdot\nabla\zeta\dd\le0
\qquad
\text{for every }0\le\zeta\in C_c^\infty(\Omega).
\]
Harmonicity on the positivity sets is established in
\cite[Theorem~1.7]{soave2024partial1}.
Thus, in any open region where two components are positive, both are
harmonic and the third vanishes. This local description does not assert
harmonicity across their free boundaries.

\begin{remark}
The endpoint regularity of the limit and the uniform estimates for
$\alpha<3/4$ are distinct conclusions. In particular,
Theorem~\ref{thm:ST} does not assert a uniform $C^{0,3/4}$ estimate
for the penalized family.
The sharp planar example also shows why, in general, a uniform Lipschitz
bound cannot be imposed on the regular approximants in
Theorem~\ref{thm:main}: an $H^1$ limit of maps with uniformly bounded
Lipschitz constants has a Lipschitz representative.
\end{remark}

The results recalled here describe minimizing configurations and their
limits. They provide the regularity context for the density problem,
but are not used in its proof. Positive penalized solutions have
$u_1^\eps u_2^\eps u_3^\eps>0$ before taking the limit, whereas
Theorem~\ref{thm:main} constructs smooth competitors satisfying the
exact constraint at every approximation step. That construction applies
to every $U\in S$, without requiring $U$ to be a minimizer or to have
regular interfaces. The global variational results in
Section~\ref{sec:gamma} require only the standing trace assumptions;
(h1), (h2), and free-boundary structure are not additional hypotheses
for those results.

\section{Admissibility and the segregation retraction}\label{sec:retraction}

Throughout the paper, we work under the following standing assumptions, which are strengthened only where it is stated explicitly
\begin{equation}\label{eq:standing}
\begin{gathered}
\Omega\subset\R^d\ \text{is a bounded Lipschitz domain},\quad d\ge1,\\
\varphi=(\varphi_1,\varphi_2,\varphi_3)\in H^{1/2}(\partial\Omega)^3,\quad
\varphi_i\ge0\ \ (i=1,2,3)\quad\text{and}\quad
\varphi_1\varphi_2\varphi_3=0\quad\text{a.e. on }\partial\Omega.
\end{gathered}
\end{equation}
All products and relations are understood in the sense of pointwise products of almost everywhere defined representatives.

\subsection{The geometric retraction}

Set $\one=(1,1,1)$ and define the retraction map $\calR : \R^3 \to \R^3$ by
\begin{equation}\label{eq:retraction}
\calR(a)=a-\min\{a_1,a_2,a_3\}\one,
\quad a=(a_1,a_2,a_3)\in\R^3.
\end{equation}
We set 
$$
\mathcal{K} = \{ a = (a_1, a_2, a_3) \in \mathbb{R}^3 \mid a_i \ge 0 \text{ for } i=1,2,3, \text{ and } a_1 a_2 a_3 = 0 \}.
$$

\begin{lemma}[The segregation retraction]\label{lem:retraction}
\begin{enumerate}
\item[\textup{(i)}] For every $a\in\R^3$, $\calR(a)\in\mathcal K$; that is,
\begin{equation}\label{eq:retraction-constraint}
\calR_i(a)\ge0\quad(i=1,2,3),
\qquad
\calR_1(a)\,\calR_2(a)\,\calR_3(a)=0 .
\end{equation}
\item[\textup{(ii)}] $\calR(a)=a$ if and only if $a\in\mathcal K$. Consequently, $\calR(\R^3)=\mathcal K$ and $\calR\circ\calR=\calR$, meaning $\calR$ is a retraction of $\R^3$ onto $\mathcal K$.
\item[\textup{(iii)}] The map $\calR$ is globally Lipschitz, satisfying
\begin{equation}\label{eq:retraction-lip}
|\calR(a)-\calR(b)|\le\sqrt3\,|a-b|
\qquad\text{for all }a,b\in\R^3.
\end{equation}
The constant $\sqrt3$ is optimal, with equality achieved for $a=(0,1,1)$ and $b=(2t,1-t,1-t)$ where $0<t\le\tfrac13$.
\item[\textup{(iv)}] For every $V\in H^1(\Omega)^3$, the composition satisfies $\calR(V)\in H^1(\Omega)^3$ with
\begin{equation}\label{eq:retraction-gradient}
|\nabla\calR(V)|\le\sqrt3\,|\nabla V|
\quad\text{a.e. in }\Omega.
\end{equation}
Consequently, $E(\calR(V))\le3\,E(V)$.
\end{enumerate}
\end{lemma}

\begin{proof}
\textit{Proof of (i) and (ii).} Let $m(a) \coloneqq \min\{a_1, a_2, a_3\}$. By definition, $\calR_i(a) = a_i - m(a) \ge 0$ for all $i \in \{1,2,3\}$. Furthermore, $a_k = m(a)$ for at least one index $k$, which implies $\calR_k(a) = 0$. This ensures $\calR_1(a)\calR_2(a)\calR_3(a) = 0$, verifying that $\calR(a) \in \mathcal{K}$. 
If $a \in \mathcal{K}$, then $a_i \ge 0$ for all $i$ and $a_k = 0$ for some $k$. This yields $m(a) = 0$, so $\calR(a) = a$. Conversely, if $\calR(a) = a$, then $m(a)\one = 0$, implying $m(a) = 0$ and thus $a \in \mathcal{K}$. This establishes that the fixed-point set of $\calR$ is precisely $\mathcal{K}$. Combined with the property $\calR(\R^3) \subset \mathcal{K}$, we deduce that $\calR(\R^3) = \mathcal{K}$ and $\calR \circ \calR = \calR$.

\medskip
\textit{Proof of (iii).} For each $k \in \{1,2,3\}$, define the closed convex cone
\[
C_k \coloneqq \bigl\{a\in\R^3:\ a_k\le a_j\ \text{for } j=1,2,3\bigr\}.
\]
Observe that $\bigcup_{k=1}^3 C_k = \R^3$. On each cone $C_k$, we have $m(a) = a_k$, and thus $\calR$ coincides with the linear operator $J_k \coloneqq I - \one e_k^{\top} \in \R^{3 \times 3}$. Note that on any intersection $C_k \cap C_{k'}$, the relation $a_k = a_{k'}$ holds, ensuring $\calR$ is well-defined and continuous across boundaries.

To determine the operator norm of $J_k$, we evaluate $|J_k \xi|^2 = \sum_{i=1}^3 (\xi_i - \xi_k)^2$ for $\xi \in \R^3$. For $k=1$, this corresponds to the quadratic form $\xi^{\top} J_1^{\top} J_1 \xi$, where
\[
J_1^{\top} J_1 = \begin{pmatrix} 2 & -1 & -1 \\ -1 & 1 & 0 \\ -1 & 0 & 1 \end{pmatrix}.
\]
The characteristic polynomial $\det(J_1^{\top} J_1 - \lambda I) = -\lambda(1-\lambda)(3-\lambda)$ yields eigenvalues $0, 1$, and $3$. Consequently, $|J_1 \xi|^2 \le 3|\xi|^2$, with equality holding if and only if $\xi$ is proportional to the eigenvector $(2,-1,-1)^{\top}$. Since $J_k = P_k J_1 P_k$, where $P_k$ is the permutation matrix exchanging coordinates $1$ and $k$, symmetry implies that
\begin{equation}\label{eq:Jk-bound}
|J_k \xi| \le \sqrt{3}\,|\xi| \qquad \text{for all } \xi \in \R^3 \text{ and } k \in \{1,2,3\},
\end{equation}
with equality achieved when $\xi$ is parallel to $3e_k - \one$.

To establish the global Lipschitz bound, consider $a,b \in \R^3$ and the line segment $\gamma(s) \coloneqq a + s(b-a)$ for $s \in [0,1]$. For each $k$, the preimage $I_k \coloneqq \{s \in [0,1] : \gamma(s) \in C_k\}$ is a closed interval (possibly degenerate or empty), as it is the preimage of a closed convex set under an affine map. There exists a finite partition $0 = s_0 < s_1 < \dots < s_N = 1$ comprising $0, 1$, and the endpoints of the non-empty intervals $I_k$. Each open subinterval $(s_{p-1}, s_p)$ is contained entirely within some $I_{k_p}$, and by continuity, so is its closure $[s_{p-1}, s_p]$. On this segment, $\calR \circ \gamma = J_{k_p} \circ \gamma$. Applying \eqref{eq:Jk-bound} yields
\[
\bigl|\calR(\gamma(s_p)) - \calR(\gamma(s_{p-1}))\bigr|
= \bigl|J_{k_p}\bigl(\gamma(s_p) - \gamma(s_{p-1})\bigr)\bigr|
\le \sqrt{3}\,(s_p - s_{p-1})|b-a|.
\]
Summing over $p = 1, \dots, N$ and applying the triangle inequality yields \eqref{eq:retraction-lip}.

To verify optimality, set $a = (0,1,1)$ and $b = (2t, 1-t, 1-t)$ for $t \in (0, 1/3]$. Since $2t \le 1-t$, both points lie in $C_1$. We compute $\calR(a) = (0,1,1)$ and $\calR(b) = (0, 1-3t, 1-3t)$, yielding
\[
\frac{|\calR(a) - \calR(b)|}{|a - b|} = \frac{3\sqrt{2}\,t}{\sqrt{6}\,t} = \sqrt{3}.
\]
Thus, the constant $\sqrt{3}$ cannot be improved.

\medskip
\textit{Proof of (iv).} Let $V = (v_1, v_2, v_3) \in H^1(\Omega)^3$. By the lattice properties of $H^1(\Omega)$ \cite[Ch.~II, Appendix~A]{KinderlehrerStampacchia1980}, $f \in H^1(\Omega)$ implies $|f| \in H^1(\Omega)$, with weak gradient $\nabla|f| = \nabla f \one_{\{f>0\}} - \nabla f \one_{\{f<0\}}$ almost everywhere. Utilizing the identity $\min\{f,g\} = \frac{1}{2}(f + g - |f - g|)$, we deduce that $m(V) \coloneqq \min\{v_1, v_2, v_3\} \in H^1(\Omega)$. Consequently, $\calR(V) = V - m(V)\one \in H^1(\Omega)^3$.

To evaluate the weak gradient, we partition $\Omega$ into the disjoint measurable sets:
\begin{align*}
A_1 &\coloneqq \{v_1 \le v_2\} \cap \{v_1 \le v_3\}, \\
A_2 &\coloneqq \{v_2 < v_1\} \cap \{v_2 \le v_3\}, \\
A_3 &\coloneqq \{v_3 < v_1\} \cap \{v_3 < v_2\}.
\end{align*}
On each $A_k$, $V(x) \in C_k$ and $m(V) = v_k$. Recall that for any $f, g \in H^1(\Omega)$, $\nabla f = \nabla g$ almost everywhere on the set $\{f=g\}$. Applying this to $f = m(V)$ and $g = v_k$ implies $\nabla m(V) = \nabla v_k$ a.e.\ on $A_k$. Thus, for each partial derivative $\partial_\ell$ ($\ell \in \{1,\dots,d\}$), we have
\[
\partial_\ell \calR(V) = \partial_\ell V - (\partial_\ell v_k)\one = J_k \partial_\ell V \quad \text{a.e. on } A_k,
\]
where $\partial_\ell V \coloneqq (\partial_\ell v_1, \partial_\ell v_2, \partial_\ell v_3)^{\top} \in \R^3$. Applying \eqref{eq:Jk-bound} to each of the $d$ columns of the Jacobian matrix $\nabla \calR(V)$ yields
\[
|\nabla \calR(V)|^2 = \sum_{\ell=1}^d |J_k \partial_\ell V|^2 \le 3 \sum_{\ell=1}^d |\partial_\ell V|^2 = 3 |\nabla V|^2 \quad \text{a.e. on } A_k.
\]
Since $\Omega = A_1 \cup A_2 \cup A_3$ up to a set of measure zero, the pointwise bound \eqref{eq:retraction-gradient} holds globally almost everywhere. Integrating over $\Omega$ immediately yields $E(\calR(V)) \le 3E(V)$.
\end{proof}

\begin{lemma}[Strong continuity of the retraction]\label{lem:strong-continuity}
If a sequence $V^n\to V$ strongly in $H^1(\Omega)^3$, then
\begin{equation}\label{eq:strong-retraction}
\calR(V^n)\to\calR(V)
\quad\text{strongly in }H^1(\Omega)^3.
\end{equation}
Furthermore, the trace operator commutes with the retraction:
\begin{equation}\label{eq:trace-commutes}
\Tr\calR(V)=\calR(\Tr V).
\end{equation}
\end{lemma}

\begin{proof}
By the  properties of Sobolev spaces (see \cite[Chapter II]{KinderlehrerStampacchia1980}), the absolute value map $f \mapsto |f|$ is strongly continuous in $H^1(\Omega)$. Through the algebraic identity
\[
\min\{f,g\}=\frac{f+g-|f-g|}{2},
\]
this continuity extends to the pointwise minimum of any finite number of $H^1$ functions. Because the retraction $\calR_i(V) = v_i - \min_j v_j$ is a linear combination of the identity and the minimum operator, the strong convergence $V^n \to V$ guarantees $\calR(V^n) \to \calR(V)$ in $H^1(\Omega)^3$.

Finally, we establish the trace identity \eqref{eq:trace-commutes} by density. Let $V \in H^1(\Omega)^3$. Since $\Omega$ is Lipschitz, there exists a sequence $\{V^n\} \subset C^\infty(\overline{\Omega})^3$ such that 
\[
V^n \to V \quad \text{strongly in } H^1(\Omega)^3.
\]
For functions smooth up to the boundary, the trace is the classical pointwise restriction. Thus, the continuous retraction commutes on $\partial\Omega$ as follows
\begin{equation}\label{eq:trace-smooth}
\Tr \calR(V^n) = \calR(\Tr V^n) \quad \text{pointwise on } \partial\Omega.
\end{equation}
The standard trace operator $\Tr : H^1(\Omega) \to L^2(\partial\Omega)$ is linear and bounded. By \eqref{eq:strong-retraction}, we know $\calR(V^n) \to \calR(V)$ in $H^1(\Omega)^3$, so applying the trace yields
\begin{equation}\label{eq:trace-limit-1}
\Tr \calR(V^n) \to \Tr \calR(V) \quad \text{strongly in } L^2(\partial\Omega)^3.
\end{equation}
Simultaneously, because $\calR:\R^3\to\R^3$ is globally Lipschitz (with constant $\sqrt{3}$), it induces a continuous Nemytskii operator on $L^2(\partial\Omega)^3$. Since $\Tr V^n \to \Tr V$ in $L^2(\partial\Omega)^3$, we obtain
\begin{equation}\label{eq:trace-limit-2}
\norm{\calR(\Tr V^n) - \calR(\Tr V)}_{L^2(\partial\Omega)^3} \le \sqrt{3} \norm{\Tr V^n - \Tr V}_{L^2(\partial\Omega)^3} \to 0.
\end{equation}
Equating the $L^2(\partial\Omega)^3$ limits of the left and right sides of \eqref{eq:trace-smooth} via \eqref{eq:trace-limit-1} and \eqref{eq:trace-limit-2}  implies
\[
\Tr \calR(V) = \calR(\Tr V). \qedhere
\]
\end{proof}

\subsection{Nonemptiness and existence of minimizers}

\begin{proposition}[Compatible traces ensure admissibility]\label{prop:nonempty}
Under \eqref{eq:standing}, $S\ne\varnothing$. Under the stronger condition \eqref{eq:lip-data}, there exists
\begin{equation}\label{eq:lip-lifting}
\Phi\in S\cap\Lip(\overline\Omega)^3.
\end{equation}
\end{proposition}

\begin{proof}
Under the standing assumption \eqref{eq:standing}, each component of the boundary data satisfies $\varphi_i \in H^{1/2}(\partial\Omega)$. For any bounded Lipschitz domain, the classical trace operator $\Tr : H^1(\Omega) \to H^{1/2}(\partial\Omega)$ is surjective. By the definition of surjectivity, every element in the target space is the trace of some function in the domain; thus, there exist scalar Sobolev functions $f_i \in H^1(\Omega)$ such that $\Tr f_i = \varphi_i$ for each $i=1,2,3$. 

Gathering these components into a vector field $F = (f_1, f_2, f_3) \in H^1(\Omega)^3$, we construct the candidate configuration by 
\[
\Phi := \calR(F).
\]
By Lemma~\ref{lem:retraction}, $\Phi \in H^1(\Omega)^3$, its components are non-negative, and they satisfy the exclusion constraint $\Phi_1\Phi_2\Phi_3 = 0$ almost everywhere in $\Omega$. Furthermore, by the trace commutativity established in Lemma~\ref{lem:strong-continuity}, we can evaluate the boundary values of this construction:
\[
\Tr \Phi = \calR(\Tr F) = \calR(\varphi).
\]
Because the prescribed boundary data satisfies $\varphi_i \ge 0$ and $\varphi_1\varphi_2\varphi_3 = 0$ almost everywhere on $\partial\Omega$, the retraction $\calR$ acts as the identity mapping on $\varphi$. Therefore, $\calR(\varphi) = \varphi$, which confirms that $\Tr \Phi = \varphi$. This verifies that $\Phi \in S$, proving the constrained class is non-empty.

If we assume the regularity condition $\varphi \in \Lip(\partial\Omega)^3$, we may extend its components to Lipschitz functions defined on all of $\R^d$ via the classical McShane-Whitney extension theorem (\cite[Theorem~3.1]{EvansGariepy2015}). Applying the Lipschitz retraction $\calR$ to the resulting vector-valued extension yields \eqref{eq:lip-lifting}.
\end{proof}

\begin{lemma}\label{lem:weakclosed}
The space $S$ is sequentially weakly closed in $H^1(\Omega)^3$.
\end{lemma}

\begin{proof}
Let $U^n\in S$ be a sequence such that $U^n\rightharpoonup U$ weakly in $H^1(\Omega)^3$. By the Rellich-Kondrachov Theorem, $H^1(\Omega)$ is compactly embedded in $L^2(\Omega)$, providing a strongly convergent subsequence in $L^2(\Omega)^3$, which we may further extract so that
\[
U^n\to U\quad\text{a.e. in }\Omega.
\]
Because both the nonnegativity $u_i^n \ge 0$ and the exclusion constraint $u_1^n u_2^n u_3^n = 0$ hold pointwise a.e., they pass to the almost everywhere limit $U$. Furthermore, the bounded linear trace operator is weakly continuous, ensuring $\Tr U=\varphi$. Thus $U\in S$.
\end{proof}

\begin{proposition}\label{prop:existence}
Under \eqref{eq:standing}, the constrained energy $E^0$ and every penalized energy $E^\eps$ ($\eps>0$) admit minimizers with finite energy.
\end{proposition}

\begin{proof}
The existence of minimizers follows from the direct method of the calculus of variations. Let $\Phi \in S$ be the admissible lifting provided by Proposition~\ref{prop:nonempty}. Because $\Phi$  satisfies the exclusion constraint, its penalization integral vanishes, making it a common finite-energy competitor:
\[
E^\eps(\Phi) = E^0(\Phi) = E(\Phi) < \infty.
\]
Consequently, the infima $m_0 \coloneqq \inf_{S} E^0$ and $m_\eps \coloneqq \inf_{\mathcal{A}} E^\eps$ are both finite and bounded below by zero.

Let $\{U^n\}_{n=1}^\infty$ be a minimizing sequence for either the constrained or penalized problem. In either case, the unpenalized Dirichlet energy is uniformly bounded, i.e. $E(U^n) \le C$. Because $\Tr U^n = \Tr \Phi = \varphi$, the difference satisfies $U^n - \Phi \in H_0^1(\Omega)^3$. Applying the Poincaré inequality to this difference yields 
\begin{equation}\label{eq:coercivity}
\begin{aligned}
\norm{U^n}_{H^1(\Omega)^3} 
&\le \norm{U^n - \Phi}_{H^1(\Omega)^3} + \norm{\Phi}_{H^1(\Omega)^3} \\
&\le C_{\Omega} \norm{\nabla U^n - \nabla \Phi}_{L^2(\Omega)^3} + \norm{\Phi}_{H^1(\Omega)^3} \\
&\le C_{\Omega} E(U^n)^{1/2} + (1+C_{\Omega})\norm{\Phi}_{H^1(\Omega)^3}.
\end{aligned}
\end{equation}
This guarantees that any minimizing sequence is uniformly bounded in $H^1(\Omega)^3$. 

By the Banach-Alaoglu and Rellich-Kondrachov compactness theorems, we extract a subsequence such that
\[
U^n \rightharpoonup U \quad \text{weakly in } H^1(\Omega)^3, 
\qquad U^n \to U \quad \text{strongly in } L^2(\Omega)^3 \text{ and a.e.\ in } \Omega.
\]
Because the Dirichlet energy $E(V) = \norm{\nabla V}_{L^2(\Omega)^3}^2$ is convex and continuous, it is weakly lower semicontinuous on $H^1(\Omega)^3$, meaning:
\[
E(U) \le \liminf_{n\to\infty} E(U^n).
\]

For the constrained problem ($E^0$), the sequence $\{U^n\} \subset S$. By Lemma~\ref{lem:weakclosed}, the set $S$ is sequentially weakly closed, ensuring the limit $U \in S$. Thus, $E^0(U) = E(U) \le \liminf E^0(U^n) = m_0$, confirming $U$ is a global minimizer.

For the penalized problem ($E^\eps$), the sequence $\{U^n\} \subset \mathcal{A}$. An analogous weak closure argument guarantees $U \in \mathcal{A}$. Because $U^n \to U$ almost everywhere, the pointwise product satisfies $(u_1^n u_2^n u_3^n)^2 \to (u_1 u_2 u_3)^2$ almost everywhere in $\Omega$. Applying Fatou's lemma to this non-negative sequence of integrands yields:
\[
\int_\Omega (u_1 u_2 u_3)^2 \dd \le \liminf_{n\to\infty} \int_\Omega (u_1^n u_2^n u_3^n)^2 \dd.
\]
Summing the lower bounds for the Dirichlet energy and the penalization integral establishes that the full functional $E^\eps$ is weakly lower semicontinuous:
\[
E^\eps(U) \le \liminf_{n\to\infty} \left( E(U^n) + \frac{1}{\eps} \int_\Omega (u_1^n u_2^n u_3^n)^2 \dd \right) = \liminf_{n\to\infty} E^\eps(U^n) = m_\eps.
\]
This proves the existence of a finite-energy minimizer for every $\eps > 0$.
\end{proof}

\begin{proposition}\label{prop:bounded-minimizers}
Assume \eqref{eq:standing} and $\varphi_i\in L^\infty(\partial\Omega)$. Every global minimizer $U$ of $E^0$ or $E^\eps$ satisfies
\begin{equation}\label{eq:maximum-bound}
0\le u_i\le M_i:=\|\varphi_i\|_{L^\infty(\partial\Omega)}
\quad\text{a.e. in }\Omega.
\end{equation}
For fixed $\eps>0$, every penalized minimizer belongs to $C^\infty(\Omega)^3$ and satisfies \eqref{eq:penalized_system} in all of $\Omega$. On each connected component of $\Omega$, each $u_i^\eps$ is either identically zero or strictly positive.
\end{proposition}

\begin{proof}
Replace $u_i$ by $\min\{u_i,M_i\}$ and leave the other components unchanged. The trace, nonnegativity, and, for $E^0$, segregation are preserved. The penalty cannot increase, while the Dirichlet energy decreases by
\[
\int_{\{u_i>M_i\}}|\nabla u_i|^2\dd
=\|\nabla(u_i-M_i)_+\|_{L^2}^2.
\]
Minimality makes this quantity zero. Since $(u_i-M_i)_+\in H_0^1(\Omega)$, Poincar\'e's inequality proves \eqref{eq:maximum-bound}.

For the penalized problem, a signed map $W\in H^1(\Omega)^3$ with trace $\varphi$ satisfies $|W|\in\mathcal A$ and
\[
E^\eps(|W|)=\int_\Omega\left(\sum_i|\nabla w_i|^2+\eps^{-1}\prod_i w_i^2\right)\dd.
\]
Thus minimality applies to this energy expression along signed compactly supported variations. Boundedness justifies differentiation along every $\zeta\in C_c^\infty(\Omega)$, giving
\begin{equation}\label{eq:weak-EL}
\int_\Omega\nabla u_i^\eps\cdot\nabla\zeta\dd
+\frac1\eps\int_\Omega u_i^\eps\prod_{j\ne i}(u_j^\eps)^2\zeta\dd=0.
\end{equation}
The right-hand side of the equation for $\Delta u_i^\eps$ is bounded. Interior elliptic estimates give $W^{2,p}_{\mathrm{loc}}$ regularity for every finite $p$, and bootstrapping gives smoothness. Finally,
\[
-\Delta u_i^\eps+c_i^\eps u_i^\eps=0,
\quad c_i^\eps=\eps^{-1}\prod_{j\ne i}(u_j^\eps)^2\ge0,
\]
so the strong maximum principle gives the alternative.
\end{proof}

\begin{remark}
For the global minimizers considered here, bounded boundary data imply (h1), and global minimality implies (h2). On a connected domain a nontrivial boundary trace gives a positive component; if a component vanishes identically, the remaining equations are harmonic. Thus Lipschitz boundary data require no separate uniform $L^\infty$ assumption when applying the local estimates recalled above.
\end{remark}

\section{Gamma convergence and penalized minimizers}\label{sec:gamma}
Throughout this section, we assume only the standing hypotheses \eqref{eq:standing}. The $\Gamma$-convergence of the family $\{E^\eps\}_{\eps>0}$ as $\eps\to0$ is evaluated in the $L^2(\Omega)^3$ topology \cite[Chapter~1]{braides2002gamma}. This topological choice is natural since by \eqref{eq:coercivity}, energy-bounded sets are uniformly bounded in $H^1(\Omega)^3$, rendering them relatively compact in $L^2(\Omega)^3$. Because $E^\eps \le E^0$ globally with equality on $S$, the recovery sequence condition is satisfied trivially by constant sequences. Thus, the core of Theorem~\ref{thm:gamma} is the liminf inequality, which relies on the sequential weak closure of the admissible class $\mathcal A$ in $H^1(\Omega)^3$ (a property established analogously to Lemma~\ref{lem:weakclosed}). Following the $\Gamma$-convergence framework, we strengthen the $L^2$ convergence of penalized almost-minimizers to strong $H^1$ convergence and demonstrate the vanishing of the penalization energy. The construction of regular recovery sequences is discussed in Section~\ref{sec:density}.

\begin{theorem}[Gamma convergence]\label{thm:gamma}
Under \eqref{eq:standing}, the penalized functional $\Gamma$-converges to the constrained functional:
\begin{equation}\label{eq:gamma}
E^\eps\xrightarrow{\Gamma}E^0
\quad\text{in }L^2(\Omega)^3
\quad\text{as }\eps\to0.
\end{equation}
Furthermore, the family $\{E^\eps\}_{\eps>0}$ is equicoercive in $L^2(\Omega)^3$.
\end{theorem}

\begin{proof}
\emph{Compactness and Equicoercivity.} Because $E(U)\le E^\eps(U)$ for all $\eps > 0$, the uniform coercivity bound \eqref{eq:coercivity} guarantees that the union of any fixed sublevel set $\bigcup_{\eps>0} \{U : E^\eps(U) \le M\}$ is uniformly bounded in $H^1(\Omega)^3$. The Rellich-Kondrachov compactness theorem therefore asserts that this set is relatively compact in $L^2(\Omega)^3$.

\emph{Liminf inequality.} Let $\eps_n\to0$ and consider any sequence $U^n\to U$ strongly in $L^2(\Omega)^3$. If the liminf of the energies is infinite, the inequality holds trivially. If finite, we extract a subsequence attaining this limit, ensuring
\[
E^{\eps_n}(U^n)\le M<\infty.
\]
Equicoercivity provides a further subsequence such that
\[
U^n\rightharpoonup U\quad\text{in }H^1(\Omega)^3,
\quad U^n\to U\quad\text{a.e. in }\Omega.
\]
Because $\mathcal A$ is weakly closed, $U \in\mathcal A$. The penalization bound implies
\begin{equation}\label{eq:overlap-vanishes}
\int_\Omega(u_1^nu_2^nu_3^n)^2\dd\le M\eps_n\to0.
\end{equation}
By Fatou's lemma on the a.e. converging integrands,
\[
\int_\Omega(u_1u_2u_3)^2\dd = 0,
\]
which places the limit inside the constraint set $U\in S$. Therefore, we may drop the penalization term to conclude
\[
E^0(U)=E(U)
\le\liminf_{n\to\infty}E(U^n)
\le\liminf_{n\to\infty}E^{\eps_n}(U^n).
\]

\emph{Limsup inequality (Recovery Sequence).} If $U\in S$, we simply choose the constant sequence $U^\eps=U$. Since $u_1 u_2 u_3 = 0$ globally, the penalization term evaluates to zero
\[
E^\eps(U^\eps)=E(U)=E^0(U)
\quad\text{for every }\eps>0.
\]
If $E^0(U)=+\infty$, the limsup inequality is automatic. This completes the proof of \eqref{eq:gamma}.
\end{proof}

\begin{theorem}[Strong convergence of minimizers]\label{thm:minimizers}
Assume \eqref{eq:standing}. Let $U^\eps$ be a global minimizer of $E^\eps$ and set
\[
m_\eps=\min_{L^2(\Omega)^3}E^\eps,
\quad
m_0=\min_{L^2(\Omega)^3}E^0.
\]
Then the energy minima converge monotonically, i.e. $m_\eps\uparrow m_0$ as $\eps\to0$. For every arbitrary null sequence $\eps_n\to0$, there exists a subsequence and a minimizer $U\in S$ of $E^0$ such that
\begin{equation}\label{eq:minimizer-strong}
U^{\eps_n}\to U\quad\text{strongly in }H^1(\Omega)^3,
\end{equation}
and the penalization energy vanishes
\begin{equation}\label{eq:penalty-vanishes}
\frac1{\eps_n}\int_\Omega
\bigl(u_1^{\eps_n}u_2^{\eps_n}u_3^{\eps_n}\bigr)^2\dd
\to0.
\end{equation}
\end{theorem}

\begin{proof}
Let $U^0\in S$ be a true minimizer of $E^0$. By definition of minimality and the non-negativity of the penalty term, we bracket the energies:
\begin{equation}\label{eq:energy-sandwich}
E(U^\eps)\le E^\eps(U^\eps)
\le E^\eps(U^0)=m_0.
\end{equation}
Equicoercivity provides a subsequence converging weakly in $H^1$ and strongly in $L^2$ to some $U$. The liminf inequality from Theorem~\ref{thm:gamma} implies $U\in S$. Therefore, by weak lower semicontinuity of the Dirichlet energy:
\[
m_0\le E(U)
\le\liminf_{n\to\infty}E(U^{\eps_n})
\le\limsup_{n\to\infty}E^{\eps_n}(U^{\eps_n})\le m_0.
\]
Because both ends of this chain are $m_0$, all intermediate quantities must coincide. Specifically:
\[
E(U^{\eps_n})\to E(U)=m_0,
\quad
E^{\eps_n}(U^{\eps_n})\to m_0.
\]
We now have weak convergence $U^{\eps_n} \rightharpoonup U$ in $H^1$ combined with the convergence of the gradient norms $\|\nabla U^{\eps_n}\|_{L^2}^2 \to \|\nabla U\|_{L^2}^2$. Weak convergence of the gradients and convergence of their $L^2$ norms imply strong gradient convergence. Together with strong $L^2$ convergence of the functions, this establishes \eqref{eq:minimizer-strong}. Subtracting the two convergent energies isolates the penalty term, proving \eqref{eq:penalty-vanishes}.

Since $E^{\eps'}\ge E^\eps$ for $0<\eps'<\eps$, the minima are nondecreasing as $\eps$ decreases. The subsequence argument applied to arbitrary null sequences uniformly identifies their limit as $m_0$.
\end{proof}

\begin{remark}\label{rem:almost-minimizers}
The same proof applies to almost-minimizers $V^\eps\in\mathcal A$ satisfying
\[
E^\eps(V^\eps)\le m_\eps+o(1).
\]
Every sequence $\eps_n\to0$ will contain a subsequence converging strongly in $H^1$ to a constrained minimizer, driving the penalization energy to zero. No uniqueness of the limiting minimizer is necessary or asserted.
\end{remark}

Let $\mathcal{M} \subset S$ denote the set of all global minimizers of the constrained functional $E^0$:
\begin{equation}\label{eq:minimizer-set}
\mathcal{M} := \bigl\{ U \in S : E^0(U) = \min_{W \in S} E^0(W) \bigr\}.
\end{equation}
For any configuration $V \in H^1(\Omega)^3$, the distance to the minimizer set $\mathcal{M}$ in the $H^1$ topology is defined in the standard way as the infimum of the norm difference
\begin{equation}\label{eq:distance-def}
\dist_{H^1(\Omega)^3}(V, \mathcal{M}) := \inf_{W \in \mathcal{M}} \norm{V - W}_{H^1(\Omega)^3}.
\end{equation}
Because $\mathcal{M}$ is strongly compact in $H^1(\Omega)^3$, this infimum is always attained by some specific minimizer $W_V \in \mathcal{M}$.

\begin{lemma}\label{cor:minimizing-set}

Assume \eqref{eq:standing}. Then $\mathcal M$ defined in \eqref{eq:minimizer-set}, is nonempty and compact in $H^1(\Omega)^3$.
Let $V^\eps\in\mathcal A$ satisfy
\[
E^\eps(V^\eps)\le m_\eps+r_\eps,
\qquad 0\le r_\eps\to0.
\]
Then,  as $\eps\to0$, we have
\begin{align}
\dist_{H^1(\Omega)^3}(V^\eps,\mathcal M)&\to0,\label{eq:distance-minimizers}\\
\|V^\eps-\calR(V^\eps)\|_{H^1(\Omega)^3}&\to0,\label{eq:strong-correction}\\
E^0(\calR(V^\eps))&\to m_0,\label{eq:corrected-minimum}\\
\eps^{-1}\int_\Omega(v_1^\eps v_2^\eps v_3^\eps)^2\dd&\to0.\label{eq:whole-penalty}
\end{align}
If $\mathcal M=\{U\}$, then $V^\eps\to U$ strongly in $H^1(\Omega)^3$.
\end{lemma}

\begin{proof}
The proof proceeds by establishing the topological compactness of the minimizer set $\mathcal{M}$, analyzing the subsequential convergence of the penalized configurations, and upgrading these limits to the full continuous family via a contradiction argument.

Let $\{U^k\} \subset \mathcal{M}$ be any sequence of global constrained minimizers, meaning $E^0(U^k) = m_0$ and $U^k \in S$ for all $k$. Because the Dirichlet energy is uniformly bounded, the Poincaré inequality implies that the sequence $\{U^k\}$ is uniformly bounded in the $H^1(\Omega)^3$ norm (as established by the coercivity bound \eqref{eq:coercivity}). By the Banach-Alaoglu theorem, we can extract a weakly convergent subsequence, $U^k \rightharpoonup U$ weakly in $H^1(\Omega)^3$. 

Because the constrained space $S$ is sequentially weakly closed (Lemma~\ref{lem:weakclosed}), the limit configuration satisfies the geometric exclusion and boundary conditions, placing $U \in S$. The weak lower semicontinuity of the Dirichlet energy implies
\[
E^0(U) \le \liminf_{k\to\infty} E^0(U^k) = m_0.
\]
Since $m_0$ defines the global minimum on $S$, we must have $E^0(U) = m_0$, confirming that the limit $U$ is a constrained minimizer ($U \in \mathcal{M}$). Furthermore, $E^0(U^k)\to E^0(U)$, so weak convergence of the gradients and convergence of their $L^2$ norms give
\[
\norm{\nabla(U^k-U)}_{L^2(\Omega)}\to0.
\]
Since $U^k-U\in H_0^1(\Omega)^3$, Poincar\'e's inequality yields
\[
\norm{U^k-U}_{H^1(\Omega)^3}
\le C_\Omega\norm{\nabla(U^k-U)}_{L^2(\Omega)}
\to0.
\]
Thus $\mathcal M$ is compact in $H^1(\Omega)^3$.

Let $V^\eps$ be a family of minimizers (or almost-minimizers) associated with the penalized functional $E^\eps$. By Remark~\ref{rem:almost-minimizers} (and Theorem~\ref{thm:minimizers}), for any null sequence $\eps_n \to0$, we can extract a further subsequence (which we do not relabel) such that $V^{\eps_n}$ converges strongly in $H^1(\Omega)^3$ to some target limit $U$. As proven previously, this limit is a constrained minimizer, so $U \in \mathcal{M}$. Simultaneously, the penalty term tends to zero
\[
\frac{1}{\eps_n} \int_\Omega (v_1^{\eps_n} v_2^{\eps_n} v_3^{\eps_n})^2 \dd \to 0.
\]

We now apply a retraction map $\calR$ to the converging sequence. By Lemma~\ref{lem:strong-continuity}, $\calR$ acts as a strongly continuous operator with respect to the $H^1$ topology. Therefore, the strong convergence $V^{\eps_n} \to U$ implies that the retracted sequence converges strongly as well
\[
\calR(V^{\eps_n}) \to \calR(U) \quad \text{strongly in } H^1(\Omega)^3.
\]
Because the target limit $U \in \mathcal{M} \subset S$ satisfies the partial segregation constraint $u_1 u_2 u_3 \equiv 0$, the retraction acts as the identity mapping on $U$, yielding $\calR(U) = U$. This forces the distance between the perturbed configuration and its correction to vanish in $H^1$ norm
\[
\norm{V^{\eps_n} - \calR(V^{\eps_n})}_{H^1} \le \norm{V^{\eps_n} - U}_{H^1} + \norm{U - \calR(V^{\eps_n})}_{H^1} \to 0.
\]
Moreover, since $U\in\mathcal M$,
\[
\dist_{H^1(\Omega)^3}(V^{\eps_n},\mathcal M)
\le\norm{V^{\eps_n}-U}_{H^1(\Omega)^3}\to0.
\]
Strong $H^1$ convergence of the retracted sequence also gives
\[
E^0(\calR(V^{\eps_n}))
=E(\calR(V^{\eps_n}))\to E(U)=m_0.
\]
Thus all four limits in \eqref{eq:distance-minimizers}--\eqref{eq:whole-penalty} hold along this extracted subsequence.

To elevate these subsequential limits to the continuous parameter family as $\eps \to0$, we invoke the Urysohn subsequence principle via a contradiction argument. Suppose the full-family limit fails for the distance metric $\dist_{{H^1(\Omega)^3}}(V^\eps, \mathcal{M})$. Then, there must exist a separation threshold $\eta > 0$ and a sequence $\eps_k \to0$ such that
\[
\dist_{{H^1(\Omega)^3}}(V^{\eps_k}, \mathcal{M}) \ge \eta > 0 \quad \text{for all } k.
\]
However, applying the same coercivity and compactness framework to this sequence $\eps_k$, we can extract a sub-subsequence $\eps_{k_j}$ such that $V^{\eps_{k_j}}$ converges strongly to some $\widetilde{U} \in \mathcal{M}$. This implies $\dist_{{H^1(\Omega)^3}}(V^{\eps_{k_j}},\mathcal M)\to0$, contradicting the lower bound $\eta$. 

The same contradiction argument applies to each of the remaining quantities in \eqref{eq:strong-correction}--\eqref{eq:whole-penalty}, i.e., every null sequence has a further subsequence along which all the stated limits hold. Therefore, all four limits hold as $\eps\to0$, without extraction. Finally, if the constrained minimizer is unique ($\mathcal{M} = \{U\}$), the vanishing distance to the compact set $\mathcal{M}$ in \eqref{eq:distance-minimizers} forces the continuous family $V^\eps$ to converge to that unique configuration $U$, completing the proof.
\end{proof}

\section{Smooth approximation inside the constrained class}\label{sec:density}

Throughout this section, the boundary data satisfy the Lipschitz condition \eqref{eq:lip-data}. The proof of Theorem~\ref{thm:main} separates the transition to Lipschitz continuity from the subsequent interior smoothing.

\subsection{Auxiliary results}

Before proving the main approximation theorem, we establish three key results. First, we utilize the geometric retraction to show that globally Lipschitz segregated configurations are strongly dense in the constrained Sobolev space. Next, to protect the boundary trace during smoothing, we construct a variable spatial radius $h_\delta$ such that
\[
0<h_\delta(x)\le\delta\,\operatorname{dist}(x,\partial\Omega).
\]
Finally, we design a truncated convolution procedure that smooths these Lipschitz intermediate maps while strictly preserving their zero level sets.

\begin{lemma}\label{lem:lip-density}
For every $U\in S$, there exists a sequence $W^n\in S\cap\Lip(\overline\Omega)^3$ such that
\begin{equation}\label{eq:lip-density}
W^n\to U\quad\text{strongly in }H^1(\Omega)^3.
\end{equation}
\end{lemma}

\begin{proof}
Let $U \in S$. By Proposition~\ref{prop:nonempty}, there exists a globally Lipschitz lifting $\Phi \in S \cap \Lip(\overline\Omega)^3$ such that $\Tr \Phi = \varphi$. 
Since $U \in S$, we have $\Tr U = \varphi$, which implies $\Tr(U - \Phi) = 0$. Thus,
\[
U - \Phi \in H^1_0(\Omega)^3.
\]
By the density of $C_c^\infty(\Omega)$ in $H^1_0(\Omega)$, there exists a sequence $\{\Psi^n\}_{n=1}^\infty \subset C_c^\infty(\Omega)^3$ satisfying
\begin{equation}\label{eq:smooth-approx}
\Psi^n \to U - \Phi \quad \text{strongly in } H^1(\Omega)^3.
\end{equation}
Define 
\[
F^n := \Phi + \Psi^n.
\]
Since $\Phi \in \Lip(\overline\Omega)^3$ and $\Psi^n \in C_c^\infty(\Omega)^3$, it follows that $F^n \in \Lip(\overline\Omega)^3$. From \eqref{eq:smooth-approx}, we obtain:
\[
F^n \to U \quad \text{strongly in } H^1(\Omega)^3.
\]
To restore the partial segregation constraint, we apply the retraction $\calR$ defined in \eqref{eq:retraction}, setting
\[
W^n := \calR(F^n).
\]
We verify the three required properties of $W^n$:

\begin{enumerate}
    \item \textbf{Regularity ($W^n \in \Lip$):} Because $\calR:\R^3\to\R^3$ is globally Lipschitz (Lemma~\ref{lem:retraction}) and $F^n \in \Lip(\overline\Omega)^3$, the composition maintains Lipschitz continuity: $W^n \in \Lip(\overline\Omega)^3$.
    
    \item \textbf{Admissibility ($W^n \in S$):} The properties of $\calR$ established in \eqref{eq:retraction-constraint} guarantee $w^n_i \ge 0$ and $w^n_1 w^n_2 w^n_3 \equiv 0$ in $\Omega$. Since $\Psi^n$ is compactly supported in $\Omega$, $\Tr \Psi^n = 0$, yielding $\Tr F^n = \Tr \Phi = \varphi$. By the trace commutativity in Lemma~\ref{lem:strong-continuity}:
    \[
    \Tr W^n = \calR(\Tr F^n) = \calR(\varphi).
    \]
    Because $\varphi_i \ge 0$ and $\varphi_1\varphi_2\varphi_3=0$, the retraction acts as the identity on the boundary: $\calR(\varphi) = \varphi$. Thus, $\Tr W^n = \varphi$, confirming $W^n \in S$.
    
    \item \textbf{Convergence ($W^n \to U$):} By Lemma~\ref{lem:strong-continuity}, the operator $\calR$ is strongly continuous in the $H^1$ topology. Passing the limit through the retraction yields
    \[
    W^n = \calR(F^n) \to \calR(U) \quad \text{strongly in } H^1(\Omega)^3.
    \]
    Since $U \in S$, its components satisfy $u_i \ge 0$ and $u_1 u_2 u_3 = 0$ almost everywhere. Therefore, the minimum component is zero almost everywhere, which forces $\calR(U) = U$. Consequently, $W^n \to U$ strongly in $H^1(\Omega)^3$. \qedhere
\end{enumerate}
\end{proof}

Any smoothing operation must remain confined to the interior of $\Omega$, with the smoothing radius decaying to zero near the boundary to avoid perturbing the prescribed trace. No regularity of the distance function $\dist(x, \partial\Omega)$ is required.

\begin{lemma}\label{lem:radius}
There exists a smooth cutoff function $\eta\in C^\infty(\R^d)$ such that
\begin{equation}\label{eq:eta}
\eta>0\ \text{in }\Omega,
\quad \eta=0\ \text{on }\R^d\setminus\Omega,
\quad \norm{\nabla\eta}_{L^\infty(\R^d)}\le1,
\end{equation}
and globally bounded by the Euclidean distance
\begin{equation}\label{eq:eta-distance}
0<\eta(x)\le\dist(x,\partial\Omega)
\quad \forall x\in\Omega.
\end{equation}
\end{lemma}

\begin{proof}
Let $\{K_j\}_{j=1}^\infty$ be a sequence of compact sets exhausting $\Omega$, such that $\bigcup_{j=1}^\infty K_j = \Omega$. For each $j$, construct a standard cutoff function $\psi_j \in C_c^\infty(\Omega)$ satisfying $0 \le \psi_j \le 1$ and $\psi_j > 0$ on $K_j$. We define the scaling coefficients and the target function by
\[
a_j \coloneqq \frac{2^{-j}}{1+\max_{|\alpha|\le j} \norm{D^\alpha\psi_j}_{L^\infty(\R^d)}},
\qquad
\eta_0(x) \coloneqq \sum_{j=1}^\infty a_j\psi_j(x).
\]
Fix an arbitrary multi-index $\beta$. For all $j \ge |\beta|$, the definition of $a_j$ gives
\[
\norm{D^\beta(a_j\psi_j)}_{L^\infty(\R^d)} 
= a_j \norm{D^\beta \psi_j}_{L^\infty(\R^d)} 
\le \frac{2^{-j} \norm{D^\beta \psi_j}_{L^\infty(\R^d)}}{1 + \norm{D^\beta \psi_j}_{L^\infty(\R^d)}} 
< 2^{-j}.
\]
Since $\sum_{j=1}^\infty 2^{-j} = 1 < \infty$, the Weierstrass M-test ensures that the series of partial derivatives $\sum_{j=1}^\infty D^\beta(a_j\psi_j)$ converges uniformly on $\R^d$. This uniform convergence justifies term-by-term differentiation to all orders, yielding $\eta_0 \in C^\infty(\R^d)$ with $D^\beta\eta_0 = \sum_{j=1}^\infty a_j D^\beta\psi_j$.

For any $x \in \Omega$, there exists $j$ such that $x \in K_j$, implying $a_j\psi_j(x) > 0$. Since all terms are nonnegative, $\eta_0(x) > 0$. For $x \in \R^d \setminus \Omega$, $\psi_j(x) = 0$ for all $j$ because $\operatorname{supp}(\psi_j) \subset \Omega$, ensuring $\eta_0 \equiv 0$ on $\R^d \setminus \Omega$. 

To evaluate the global Lipschitz constant, consider the first-order partial derivatives ($|\beta|=1$). For each coordinate direction $i \in \{1,\dots,d\}$,
\[
\norm{\partial_{x_i}\eta_0}_{L^\infty(\R^d)} \le \sum_{j=1}^\infty a_j \norm{\partial_{x_i}\psi_j}_{L^\infty(\R^d)} \le \sum_{j=1}^\infty 2^{-j} = 1.
\]
The Euclidean norm of the gradient is therefore uniformly bounded
\[
\norm{\nabla\eta_0}_{L^\infty(\R^d)} = \sup_{x \in \R^d} \left( \sum_{i=1}^d |\partial_{x_i}\eta_0(x)|^2 \right)^{1/2} \le \left( \sum_{i=1}^d 1^2 \right)^{1/2} = \sqrt{d}.
\]
Defining the normalized function $\eta \coloneqq \eta_0 / \sqrt{d}$ immediately scales the gradient bound to $\norm{\nabla\eta}_{L^\infty(\R^d)} \le 1$, fulfilling \eqref{eq:eta}.

Finally, to establish \eqref{eq:eta-distance}, let $x \in \Omega$. Because $\partial\Omega$ is closed, there exists a boundary point $y \in \partial\Omega$ attaining the distance infimum, satisfying $|x-y| = \dist(x,\partial\Omega)$. By construction, $\eta(y)=0$. Since $\eta \in C^\infty(\R^d)$ and possesses a bounded global gradient, the mean value theorem applied along the line segment $[x,y] \subset \R^d$ implies
\[
\eta(x) = |\eta(x) - \eta(y)| \le \left( \sup_{z \in [x,y]} |\nabla\eta(z)| \right) |x-y| \le \norm{\nabla\eta}_{L^\infty(\R^d)} |x-y| \le \dist(x,\partial\Omega).
\]
This confirms $0 < \eta(x) \le \dist(x,\partial\Omega)$ for all $x \in \Omega$, completing the proof.
\end{proof}

To prevent mollification from destroying the segregation constraint, we will truncate the mollified values where they drift too far from the original zeroes. Fix a nonnegative standard mollifier $\rho\in C_c^\infty(B_1)$ with $\int_{B_1}\rho=1$. Choose a monotonic, smooth threshold function $\theta:\R\to[0,1]$ satisfying
\[
\theta(s)=0\quad(s\le1),
\quad \theta(s)=1\quad(s\ge2),
\]
and define the active truncation map $\chi(s)=s\theta(s)$. It is not hard to check that for $s\ge0$, this map fulfills the following bounds
\begin{equation}\label{eq:chi}
\begin{gathered}
0\le\chi(s)\le s,\quad
\chi(s)=0\ (s\le1),\quad
\chi(s)=s\ (s\ge2),\\
\sup_{s\ge0}\bigl(|\chi'(s)|+|\chi(s)-s\chi'(s)|\bigr)<\infty.
\end{gathered}
\end{equation}

\begin{lemma}\label{lem:smoothing}
Let $W=(w_1,w_2,w_3)\in S\cap\Lip(\overline\Omega)^3$. There exists a family
\[
Z^\delta\in S\cap\Lip(\overline\Omega)^3\cap C^\infty(\Omega)^3,
\quad 0<\delta<\tfrac12,
\]
such that the positivity sets satisfy
\begin{equation}\label{eq:zero-preservation}
\{z_i^\delta>0\}\subseteq\{w_i>0\}
\quad\text{in }\Omega,\quad i=1,2,3,
\end{equation}
and
\begin{equation}\label{eq:smoothing-convergence}
Z^\delta\to W\quad\text{uniformly on }\overline\Omega
\quad\text{and strongly in }H^1(\Omega)^3.
\end{equation}
Furthermore, for this fixed $W$, the global Lipschitz constants of $Z^\delta$ remain bounded independently of the small parameter $\delta$.
\end{lemma}

\begin{proof}
Let $L=1+\max_i\Lip(w_i)$ denote the upper bound on the Lipschitz constants, and utilize $\eta$ from Lemma~\ref{lem:radius}. We define the spatial radius and local truncation threshold:
\begin{equation}\label{eq:variable-radius}
h_\delta(x)=\delta\eta(x),
\quad t_\delta(x)=Lh_\delta(x).
\end{equation}
For $x\in\Omega$, we execute a variable-radius convolution
\begin{equation}\label{eq:variable-mollification}
v_i^\delta(x)=\int_{B_1}\rho(y)
w_i\bigl(x-h_\delta(x)y\bigr)\,dy.
\end{equation}
Because $h_\delta(x)<\dist(x,\partial\Omega)$, the argument $x-h_\delta(x)y$ lies in $\Omega$, so $v_i^\delta\ge0$. Interior smoothness follows by differentiating the kernel in
\[
v_i^\delta(x)=\int_\Omega h_\delta(x)^{-d}
\rho\!\left(\frac{x-z}{h_\delta(x)}\right)w_i(z)\,dz,
\]
since $h_\delta$ is smooth and positive locally in $\Omega$. The Lipschitz bound gives
\begin{equation}\label{eq:mollification-error}
|v_i^\delta(x)-w_i(x)|
\le Lh_\delta(x)\int_{B_1}|y|\rho(y)\,dy
\le t_\delta(x).
\end{equation}
We now define the truncated, smooth configuration as follows
\begin{equation}\label{eq:smooth-truncation}
z_i^\delta(x):=t_\delta(x)
\chi\!\left(\frac{v_i^\delta(x)}{t_\delta(x)}\right),
\quad x\in\Omega.
\end{equation}

We show this truncation forces preservation of the zeroes. Suppose $w_i(x)=0$. Then \eqref{eq:mollification-error} gives $v_i^\delta(x)\le t_\delta(x)$. Passing this into the truncation mapping gives a ratio $v_i^\delta(x)/t_\delta(x) \le 1$. Because $\chi(s)=0$ for $s \le 1$, we have $z_i^\delta(x)=0$. This verifies \eqref{eq:zero-preservation}, meaning $\{w_i=0\}\subseteq\{z_i^\delta=0\}$. The inclusion may be strict, as new zeroes may be created near the free boundary. 

Since $W \in S$, its components satisfy $w_1 w_2 w_3 = 0$ pointwise. Therefore, at every spatial point $x \in \Omega$, at least one index $i$ must satisfy $w_i(x) = 0$. By the inclusion verified above, the corresponding $z_i^\delta(x) = 0$, enforcing partial segregation.

To establish the uniform error bound, we observe that the definition of $\theta$ implies 
$$
|\theta(s) - 1| \le 2/s\;\; \text{for all}\;\; s > 0.
$$
Because $\chi(s) = s\theta(s)$, substituting $s = q/t$ for any $q \ge 0$ and $t > 0$ yields the pointwise inequality
$$\left|t\chi(q/t)-q\right|\le 2t.$$
By substituting $q=v_i^\delta(x)$ and $t=t_\delta(x)$ into this bound, and recalling the definition of $z_i^\delta(x)$ from \eqref{eq:smooth-truncation}, we obtain
\[
|z_i^\delta(x) - v_i^\delta(x)| \le 2t_\delta(x).
\]
To bound the total deviation from the exact configuration $W$, we apply the triangle inequality 
\[
|z_i^\delta(x) - w_i(x)| 
\le |z_i^\delta(x) - v_i^\delta(x)| + |v_i^\delta(x) - w_i(x)| 
\le 2t_\delta(x) + t_\delta(x) = 3t_\delta(x).
\]
Since $t_\delta(x) = L\delta\eta(x)$ by \eqref{eq:variable-radius}, the pointwise component error is uniformly bounded by $3L\delta\norm{\eta}_{L^\infty}$. Summing the squared components for the vector-valued difference in $\R^3$ introduces a geometric factor of $\sqrt{3}$, yielding the global uniform estimate
\begin{equation}\label{eq:uniform-error}
\norm{Z^\delta - W}_{L^\infty(\Omega)^3} 
\le 3\sqrt{3} L \delta \norm{\eta}_{L^\infty(\R^d)}.
\end{equation}

Note that, as $x$ approaches the domain boundary $\partial\Omega$, the distance function forces $\eta(x) \to 0$, driving $t_\delta(x)\to0$. Extending $Z^\delta$ by the trace $W=\varphi$ on $\partial\Omega$ is therefore a continuous extension that recovers the boundary values without introducing discontinuities.

 Since $t_\delta>0$ locally in $\Omega$, \eqref{eq:smooth-truncation} is $C^\infty$ in $\Omega$. To confirm global Lipschitz bounds, we use \eqref{eq:variable-mollification} to obtain
\[
|v_i^\delta(x)-v_i^\delta(x')|
\le L\bigl(|x-x'|+|h_\delta(x)-h_\delta(x')|\bigr)
\le L(1+\delta)|x-x'|.
\]
Consider the bivariate truncation function defined on the domain $[0,\infty)^2$
\[
F(q,t) = 
\begin{cases} 
t\chi(q/t), & t>0, \\ 
q, & t=0,
\end{cases}
\qquad q \ge 0.
\]
In the interior of this domain ($t>0$), the partial derivatives evaluate to
\[
\partial_q F(q,t) = \chi'(q/t), \qquad \partial_t F(q,t) = \chi(q/t) - (q/t)\chi'(q/t).
\]
By the properties of $\chi$ established in \eqref{eq:chi}, both partial derivatives are globally bounded. This guarantees that $F$ is uniformly Lipschitz on $(0,\infty)^2$. Furthermore, the uniform error bound $|F(q,t)-q|\le 2t$ ensures that this Lipschitz continuity extends seamlessly to the boundary $t=0$.

Because the smooth truncation is formed by the composition $z_i^\delta(x) = F(v_i^\delta(x), t_\delta(x))$, we apply the chain rule for Lipschitz functions. We previously established $\Lip(v_i^\delta) \le L(1+\delta)$. Since $t_\delta(x) = L\delta\eta(x)$ with $\norm{\nabla\eta}_{L^\infty} \le 1$, we have $\Lip(t_\delta) \le L\delta$. Combining these yields
\begin{align*}
|z_i^\delta(x) - z_i^\delta(x')| 
&\le \Lip(F) \bigl( |v_i^\delta(x) - v_i^\delta(x')| + |t_\delta(x) - t_\delta(x')| \bigr) \\
&\le \Lip(F) \bigl( L(1+\delta) + L\delta \bigr) |x - x'|.
\end{align*}
Since $\delta < 1/2$, the terms depending on $\delta$ are uniformly bounded. Evaluating this across all components gives the global bound
\begin{equation}\label{eq:smooth-lip}
\Lip(Z^\delta)\le C_\chi L,
\end{equation}
where $C_\chi$ is a universal constant depending only on the chosen threshold profile $\chi$, completely independent of $\delta$ and the specific configuration $W$.

 Extend each $w_i$ to a globally Lipschitz function on $\R^d$. We calculate the weak spatial gradient of \eqref{eq:variable-mollification} via the chain rule
\begin{equation}\label{eq:mollified-gradient}
\nabla v_i^\delta(x)
=\int_{B_1}\rho(y)\nabla w_i(x-h_\delta(x)y)\,dy
-\nabla h_\delta(x)\int_{B_1}\rho(y)
\bigl[\nabla w_i(x-h_\delta(x)y)\cdot y\bigr]\,dy.
\end{equation}
To justify differentiating under the integral sign and using the chain rule, we confirm the underlying coordinate map is a homeomorphism. Because $\Lip(h_\delta)\le\delta$, an application of the triangle inequality gives, for any $|y|\le1$ and $x,x'\in\R^d$ we have
\begin{equation}\label{eq:bilipschitz-shift}
\begin{aligned}
(1-\delta)|x-x'|
\le\bigl|(x-h_\delta(x)y)-(x'-h_\delta(x')y)\bigr|
\le(1+\delta)|x-x'|.
\end{aligned}
\end{equation}
For any $\delta<1/2$, the lower bound implies the mapping is bi-Lipschitz onto its image, validating the chain rule applications.

To verify almost everywhere convergence, we recall that $\nabla w_i$ is locally integrable, meaning almost every $x \in \Omega$ is a Lebesgue point. Fix such an $x$ and denote the local smoothing radius by $h \coloneqq h_\delta(x) > 0$. Because the mollifier has unit mass ($\int_{B_1} \rho(y)\,dy = 1$), we can bring the constant $\nabla w_i(x)$ inside the integral. Applying the change of variables $z = x - hy$ (which yields $dy = h^{-d}\,dz$), we bound the difference
\begin{equation}\label{eq:lebesgue-mollification}
\begin{aligned}
\left|\int_{B_1}\rho(y)\nabla w_i(x-hy)\,dy-\nabla w_i(x)\right|
&\le \int_{B_1}\rho(y)\left|\nabla w_i(x-hy)-\nabla w_i(x)\right|\,dy \\
&\le \frac{\norm{\rho}_{L^\infty(\R^d)}}{h^d} \int_{B_h(x)}|\nabla w_i(z)-\nabla w_i(x)|\,dz.
\end{aligned}
\end{equation}
Let $\omega_d$ denote the volume of the unit ball in $\R^d$. The final term in \eqref{eq:lebesgue-mollification} is exactly $\norm{\rho}_{L^\infty(\R^d)} \omega_d$ multiplied by the integral average $\frac{1}{|B_h(x)|} \int_{B_h(x)} |\nabla w_i(z) - \nabla w_i(x)|\,dz$. As $\delta \to0$, the local radius $h = \delta\eta(x)$ approaches zero. By the Lebesgue Differentiation Theorem, this integral average vanishes.

 Furthermore, the second term in \eqref{eq:mollified-gradient} vanishes
\[
\left|\nabla h_\delta(x)\int_{B_1}\rho(y)
\bigl[\nabla w_i(x-h_\delta(x)y)\cdot y\bigr]\,dy\right|
\le L\delta.
\]
Consequently, the weak gradient converges pointwise
\begin{equation}\label{eq:gradient-ae}
\nabla v_i^\delta\to\nabla w_i\quad\text{a.e.},
\quad |\nabla v_i^\delta|\le L(1+\delta).
\end{equation}
To determine the gradient of the truncated function, we apply the chain rule to the ratio $s_i^\delta=v_i^\delta/t_\delta.$ We get
\begin{equation}\label{eq:truncated-gradient}
\nabla z_i^\delta
=\chi'(s_i^\delta)\nabla v_i^\delta
+\bigl[\chi(s_i^\delta)-s_i^\delta\chi'(s_i^\delta)\bigr]
\nabla t_\delta.
\end{equation}
We evaluate this limit in two distinct spatial regimes. First, at a point where $w_i(x)>0$, the ratio $s_i^\delta(x)$ explodes to $\infty$ because the denominator $t_\delta \to 0$ while the numerator converges to $w_i(x) > 0$. For large arguments ($s \ge 2$), $\chi(s)=s$, forcing $\chi'=1$ and $(\chi - s\chi') = 0$. Thus, the coefficients eventually become $1$ and $0$, returning $\nabla z_i^\delta \to \nabla w_i$.

Second, consider the zero set $N_i \coloneqq \{x \in \Omega : w_i(x) = 0\}$. By standard Sobolev properties, $\nabla w_i = 0$ a.e.\ on $N_i$. Consequently, \eqref{eq:gradient-ae} implies $\lim_{\delta \to 0} \nabla v_i^\delta(x) = 0$ a.e.\ on $N_i$. Simultaneously, $\norm{\nabla t_\delta}_{L^\infty} \le L\delta \to 0$. Invoking the global bounds on the truncation multipliers from \eqref{eq:chi}, taking the limit in \eqref{eq:truncated-gradient} yields
\[
\limsup_{\delta \to 0} |\nabla z_i^\delta(x)| 
\le \norm{\chi'}_\infty \lim_{\delta \to 0} |\nabla v_i^\delta(x)| 
+ \norm{\chi(s)-s\chi'(s)}_\infty \lim_{\delta \to 0} |\nabla t_\delta(x)| 
= 0.
\]
Hence, $\nabla z_i^\delta \to 0 = \nabla w_i$ a.e.\ on $N_i$.

Since $\Omega = \{w_i > 0\} \cup N_i$, we conclude that $\nabla Z^\delta \to \nabla W$ a.e.\ in $\Omega$. Furthermore, \eqref{eq:truncated-gradient} guarantees the uniform pointwise bound $\norm{\nabla Z^\delta}_{L^\infty(\Omega)^3} \le C_\chi L$. Because $|\Omega| < \infty$, this constant majorant belongs to $L^2(\Omega)$. Applying the Lebesgue Dominated Convergence Theorem gives
\[
\lim_{\delta \to 0} \norm{\nabla Z^\delta - \nabla W}_{L^2(\Omega)^3}^2 
= \int_\Omega \lim_{\delta \to 0} |\nabla Z^\delta(x) - \nabla W(x)|^2 \dd 
= 0.
\]
Coupled with the uniform estimate $\norm{Z^\delta - W}_{L^\infty(\Omega)^3} \to 0$ established in \eqref{eq:uniform-error}, this proves the strong $H^1(\Omega)^3$ convergence asserted in \eqref{eq:smoothing-convergence}.

\end{proof}

\subsection{Proof of the main  theorem}

\begin{proof}[Proof of Theorem~\ref{thm:main}]
Let $U \in S$ be a given constrained configuration. The construction of the strongly convergent smooth sequence proceeds via a diagonal extraction argument, bridging the Lipschitz approximation with the local smoothing.

By Lemma~\ref{lem:lip-density}, there exists an intermediate sequence 
\[
W^n \in S\cap\Lip(\overline\Omega)^3, \quad n \in \mathbb{N},
\]
such that $W^n$ keeps the exact boundary trace ($\Tr W^n = \varphi$) and converges to the target configuration:
\begin{equation}\label{eq:proof-Wn-limit}
\lim_{n\to\infty} \norm{W^n-U}_{H^1(\Omega)^3} = 0.
\end{equation}

Fix an integer $n \ge 1$. We apply the zero-preserving smoothing operator from Lemma~\ref{lem:smoothing} exclusively to the fixed Lipschitz map $W^n$. This generates a continuous family of parameterized functions $Z^\delta[W^n]$ for $0 < \delta < 1/2$. By \eqref{eq:smoothing-convergence} in Lemma~\ref{lem:smoothing}, we know that for this fixed $n$ the following limit holds
\[
\lim_{\delta \to0} \norm{Z^\delta[W^n] - W^n}_{H^1(\Omega)^3} = 0.
\]
Because this limit evaluates to zero, we can select a spatial parameter $\delta_n$ satisfying $0<\delta_n<\min\{1/2,1/n\}$ such that the $H^1$ error is bounded by $1/n,$ i.e.
\begin{equation}\label{eq:proof-delta-choice}
\norm{Z^{\delta_n}[W^n]-W^n}_{H^1(\Omega)^3} < \frac1n.
\end{equation}

We finalize the diagonal sequence by setting $U^n := Z^{\delta_n}[W^n]$. By Lemma~\ref{lem:smoothing}, each member $U^n$ fits into the regularity class in \eqref{eq:main-regularity}
\[
U^n \in \Lip(\overline\Omega)^3 \cap C^\infty(\Omega)^3.
\]
Furthermore, because the truncation in $Z^{\delta_n}$ preserves existing zeros (i.e., $\{w_i^n = 0\} \subseteq \{u_i^n = 0\}$), and because $W^n \in S$ guarantees $w_1^n w_2^n w_3^n = 0$ everywhere, it follows that $u_1^n u_2^n u_3^n = 0$ in $\Omega$. Coupled with the preservation of the boundary trace as the smoothing radius vanishes at $\partial\Omega$, we establish that $U^n \in S$.

To verify the strong convergence of the final sequence $U^n \to U$, we apply the standard triangle inequality in the $H^1(\Omega)^3$ norm, utilizing the bound \eqref{eq:proof-delta-choice}. We get
\[
\norm{U^n-U}_{H^1} 
\le \norm{Z^{\delta_n}[W^n] - W^n}_{H^1} + \norm{W^n-U}_{H^1} 
< \frac1n + \norm{W^n-U}_{H^1}.
\]
Taking the limit as $n \to \infty$, the term $1/n \to 0$ and, by \eqref{eq:proof-Wn-limit}, the second term also vanishes. Thus, $U^n \to U$ strongly in $H^1(\Omega)^3$, proving \eqref{eq:main-convergence}.

Let $\eps_n \to0$ be any arbitrary null sequence of penalization parameters. Because we have proven that $U^n \in S$ for all $n$, the geometric exclusion constraint is satisfied everywhere, which forces the integration of the penalty term to evaluate to zero. Therefore, the penalized energy reduces to the unconstrained Dirichlet energy
\[
E^{\eps_n}(U^n) = E^0(U^n) = E(U^n) = \int_\Omega \sum_{i=1}^3 |\nabla u_i^n|^2 \dd.
\]
Finally, to bound the difference in the Dirichlet energies, we use the difference of squares for the $L^2$ norm. Let $A = \nabla U^n$ and $B = \nabla U$. Using the reverse triangle inequality
$$
\abs{\norm{A}_{L^2}^2 - \norm{B}_{L^2}^2} \le \norm{A-B}_{L^2}(\norm{A}_{L^2} + \norm{B}_{L^2}),
$$
and applying this to our gradients yields
\begin{align*}
\abs{E(U^n)-E(U)} 
= \abs{ \norm{\nabla U^n}_{L^2}^2 - \norm{\nabla U}_{L^2}^2 }
\le \norm{\nabla U^n-\nabla U}_{L^2} \bigl(\norm{\nabla U^n}_{L^2}+\norm{\nabla U}_{L^2}\bigr).
\end{align*}
Because $U^n \to U$ strongly in $H^1$, the gradient difference $\norm{\nabla U^n-\nabla U}_{L^2} \to 0$, and the sequence of norms $\norm{\nabla U^n}_{L^2}$ is uniformly bounded. Thus, the product converges to zero, which proves \eqref{eq:main-energy} and completes the theorem.
\end{proof}

\begin{corollary}\label{cor:no-gap}
Under \eqref{eq:lip-data}, no energy gap exists between the Sobolev and smooth subclasses, i.e, the following holds
\begin{equation}\label{eq:no-gap}
\min_{U\in S}E(U)
=\inf_{U\in S\cap\Lip(\overline\Omega)^3}E(U)
=\inf_{U\in S\cap\Lip(\overline\Omega)^3\cap C^\infty(\Omega)^3}E(U).
\end{equation}
No assertion is made concerning the actual attainment of the infima over the two regular, dense subclasses.
\end{corollary}

\begin{proof}
Let us denote the three infima by
\[
m_0 \coloneqq \min_{U\in S}E(U),
\quad
m_{\text{Lip}} \coloneqq \inf_{U\in S\cap\Lip(\overline\Omega)^3}E(U),
\quad
m_\infty \coloneqq \inf_{U\in S\cap\Lip(\overline\Omega)^3\cap C^\infty(\Omega)^3}E(U).
\]
Because taking an infimum over a strictly smaller set can only increase its value, the natural subset inclusions
\[
S\cap\Lip(\overline\Omega)^3\cap C^\infty(\Omega)^3 
\subset S\cap\Lip(\overline\Omega)^3 
\subset S
\]
immediately lead to 
\begin{equation}\label{eq:proof-gap-1}
m_0 \le m_{\text{Lip}} \le m_\infty.
\end{equation}

To establish the reverse inequality, let $U^* \in S$ be a global constrained minimizer, so that $E(U^*) = m_0$. By Theorem~\ref{thm:main}, we can apply the smooth exact-segregation approximation to $U^*$. This guarantees the existence of a sequence
\[
U^n \in S\cap\Lip(\overline\Omega)^3\cap C^\infty(\Omega)^3
\]
such that $U^n \to U^*$ strongly in $H^1(\Omega)^3$. 

The Dirichlet energy $E(V) = \norm{\nabla V}_{L^2(\Omega)^3}^2$ is continuous with respect to strong $H^1$ convergence. Therefore, 
\[
\lim_{n\to\infty} E(U^n) = E(U^*) = m_0.
\]
Since each $U^n$ belongs to the most restricted functional class, the definition of the infimum requires that $m_\infty \le E(U^n)$ for every $n$. Taking the limit as $n \to \infty$ yields
\begin{equation}\label{eq:proof-gap-2}
m_\infty \le m_0.
\end{equation}
Combining \eqref{eq:proof-gap-1} and \eqref{eq:proof-gap-2} forces $m_0 = m_{\text{Lip}} = m_\infty$, confirming \eqref{eq:no-gap}.
\end{proof}

\begin{remark}\label{rem:boundary}
Compatible $H^{1/2}$ boundary data suffice for nonemptiness, global existence, and the $\Gamma$-convergence limit. The stronger Lipschitz boundary data condition \eqref{eq:lip-data} serves to accommodate the requirement that approximants be globally Lipschitz and preserve boundary values. This prerequisite is unavoidable: the boundary restriction of any globally Lipschitz function must be Lipschitz.
\end{remark}

\begin{remark}\label{rem:scale}
Theorem~\ref{thm:main} yields a family indexed by the parameter $\eps\to0$ if we re-parameterize via $U^\eps=U^{\lceil1/\eps\rceil}$. It does not assign a spatial smoothing radius to $\eps$, nor does it restrict the growth of $\Lip(U^\eps)$, or preserve the topology of the free boundary limit $U$. The preservation of existing zeros in Lemma~\ref{lem:smoothing} acts relatively upon the intermediate Lipschitz map $W^n$, not directly upon the abstract Sobolev minimizer $U$.
\end{remark}

\section{Quantitative correction of penalized configurations}\label{sec:correction}

The retraction mapping defined in \eqref{eq:retraction}  is of independent interest beyond the variational limit. It allows projecting any finite-energy penalized configuration onto the set of segregated functions without touching the boundary trace.

\begin{proposition}\label{prop:correction}
Assume \eqref{eq:standing}, and let $V\in\mathcal A$ possess finite perturbed energy $E^\eps(V)<\infty$. The retraction $\calR(V)\in S$ implies the global distance bound
\begin{equation}\label{eq:correction-distance}
\norm{V-\calR(V)}_{L^2(\Omega)^3}
\le\sqrt3\,|\Omega|^{1/3}
\bigl(\eps E^\eps(V)\bigr)^{1/6}.
\end{equation}
Furthermore, the corrected Dirichlet energy satisfies the multiplicative bound
\begin{equation}\label{eq:correction-energy}
E^0(\calR(V))\le3E(V)\le3E^\eps(V).
\end{equation}
\end{proposition}

\begin{proof}
Adherence to the segregation constraint follows from Lemma~\ref{lem:retraction}. By Lemma~\ref{lem:strong-continuity}, the trace remains unperturbed:
\[
\Tr\calR(V)=\calR(\varphi)=\varphi.
\]
To evaluate the distance, we set $m(x)=\min\{v_1(x),v_2(x),v_3(x)\}$. Because every component $v_i \ge 0$, we establish $v_i \ge m \ge 0$. Multiplying the three inequalities provides
\[
0\le m^3\le v_1v_2v_3,
\quad V-\calR(V)=m\one.
\]
By substituting the bound into the definition of the penalization energy, we extract an explicit upper bound on the $L^6$ integral
\begin{equation}\label{eq:m-six}
\int_\Omega m^6\dd
\le\int_\Omega(v_1v_2v_3)^2\dd
\le\eps E^\eps(V).
\end{equation}
To bound the $L^2$ norm by the $L^6$ norm, we apply Hölder's inequality on the Lebesgue measure
\[
\norm{V-\calR(V)}_{L^2(\Omega)^3}
=\sqrt3\norm{m}_{L^2(\Omega)}
\le\sqrt3\,|\Omega|^{1/3}\norm{m}_{L^6(\Omega)},
\]
which proves the fractional estimate \eqref{eq:correction-distance}. The energy inequality \eqref{eq:correction-energy} follows from the bounds established in \eqref{eq:retraction-gradient}.
\end{proof}

\begin{proposition}\label{prop:exact-distance}
Under \eqref{eq:standing}, we have
\begin{equation}\label{eq:L2-closure}
\overline S_{\,L^2(\Omega)^3}
=\{Y\in L^2(\Omega)^3:y_i\ge0,\ y_1y_2y_3=0\ \text{a.e.}\}.
\end{equation}
For every nonnegative $V\in L^2(\Omega)^3$, with $m=\min_i v_i$,
\begin{equation}\label{eq:exact-distance}
\dist_{L^2(\Omega)^3}(V,S)=\|m\|_{L^2(\Omega)}
=\frac1{\sqrt3}\|V-\calR(V)\|_{L^2(\Omega)^3}.
\end{equation}
In particular, if $V\in\mathcal A$ and $E^\eps(V)<\infty$,
\begin{equation}\label{eq:sharp-distance-bound}
\dist_{L^2(\Omega)^3}(V,S)
\le |\Omega|^{1/3}\left(\int_\Omega(v_1v_2v_3)^2\dd\right)^{1/6}
\le |\Omega|^{1/3}(\eps E^\eps(V))^{1/6}.
\end{equation}
Under \eqref{eq:lip-data}, the distance is unchanged if $S$ is replaced by $S\cap\Lip(\overline\Omega)^3\cap C^\infty(\Omega)^3$.
\end{proposition}

\begin{proof}
An $L^2$ limit of maps in $S$ has an almost everywhere convergent subsequence, so it is nonnegative and segregated. Conversely, let $Y$ belong to the right-hand side of \eqref{eq:L2-closure}, and fix $\Phi\in S$. Choose
\[
\Psi^n\in C_c^\infty(\Omega)^3,
\quad \Phi+\Psi^n\to Y\quad\text{in }L^2(\Omega)^3.
\]
Then $\calR(\Phi+\Psi^n)\in S$ and, by \eqref{eq:retraction-lip},
\[
\|\calR(\Phi+\Psi^n)-Y\|_{L^2}
=\|\calR(\Phi+\Psi^n)-\calR(Y)\|_{L^2}
\le\sqrt3\|\Phi+\Psi^n-Y\|_{L^2}\to0.
\]
This proves \eqref{eq:L2-closure}.

For every $W\in S$ and almost every $x\in\Omega$, there exists an index $i=i(x)$ such that $w_i(x)=0$. Therefore,
\[
|V(x)-W(x)|^2\ge v_i(x)^2\ge m(x)^2.
\]
For the reverse inequality, choose the measurable index
\[
\kappa(x)=\min\{i:v_i(x)=m(x)\},
\quad P_i(x)=v_i(x)\mathbf1_{\{i\ne\kappa(x)\}}.
\]
Then $P$ belongs to the closure in \eqref{eq:L2-closure} and
\[
|V-P|=m\quad\text{a.e.},\quad
\dist_{L^2}(V,S)\le\|V-P\|_{L^2}=\|m\|_{L^2}.
\]
Since $V-\calR(V)=m\one$, this proves \eqref{eq:exact-distance}. The bounds in \eqref{eq:sharp-distance-bound} follow from $m^6\le(v_1v_2v_3)^2$ and H\"older's inequality. The final assertion follows from Theorem~\ref{thm:main}.
\end{proof}

\begin{remark}\label{rem:nonattainment}
Let $\Omega=(-1,1)$ and $V(x)=(1+x,1-x,2)$, with $\varphi=V|_{\partial\Omega}$. These boundary data are compatible. Equality in \eqref{eq:exact-distance} would require a nearest competitor to equal
\[
P(x)=\begin{cases}
(0,1-x,2),&-1<x<0,\\
(1+x,0,2),&0<x<1.
\end{cases}
\]
While $P$ fulfills the pointwise segregation constraint, it induces a jump discontinuity
\[
\lim_{x \uparrow 0} P(x) = (0,1,2) 
\quad \text{whereas} \quad 
\lim_{x \downarrow 0} P(x) = (1,0,2).
\]
Thus, $P \notin H^1(-1,1)^3$, proving that the distance infimum is not attained within the admissible class $S$. 

This highlights the utility of the retraction $\calR(V) = V - m\one$. Because the scalar function $m(x)$ is continuous globally, $\calR(V)$ belongs to $H^1(-1,1)^3$ and provides a valid, bounded-energy competitor, whereas the strict $L^2$-nearest map fails to do so. 

Furthermore, this topological obstruction forces an energy blow-up for any sequence attempting to reach the distance infimum. Suppose $W^n \in S$ is a distance-minimizing sequence such that $\norm{V - W^n}_{L^2} \to \norm{m}_{L^2}$. This convergence in distance uniquely forces $W^n \to P$ strongly in $L^2(-1,1)^3$. If the sequence possessed uniformly bounded Dirichlet energy ($E(W^n) \le C$), the Poincaré inequality and weak sequential compactness would yield a subsequence converging weakly in $H^1$ to some $W \in S$. Weak $H^1$ convergence implies strong $L^2$ convergence, implying $W = P$ almost everywhere. This would give  $P \in H^1(-1,1)^3$, a contradiction. Thus, every distance-minimizing sequence must satisfy $E(W^n) \to \infty$.
\end{remark}

\begin{proposition}[Sharpness of the exponent]\label{prop:sharpness}
Let the boundary trace be trivial, $\varphi\equiv0$. For any nonzero, nonnegative test function $f\in C_c^\infty(\Omega)$, the parameterized family
\begin{equation}\label{eq:sharpness-family}
V^\eps=\eps^{1/6}(f,f,f),\quad 0<\eps\le1,
\end{equation}
demonstrates bounded perturbed energy and satisfies:
\begin{equation}\label{eq:sharpness-distance}
\dist_{L^2(\Omega)^3}(V^\eps,S)
=\eps^{1/6}\norm{f}_{L^2(\Omega)}.
\end{equation}
Consequently, it is impossible to derive an improved estimate of the form
\[
\dist_{L^2(\Omega)^3}(V^\eps,S)\le C\eps^\beta,
\quad \beta>\frac16,
\]
with a constant $C$ independent of $\eps$ across all bounded-energy families.
\end{proposition}

\begin{proof}
The sequence of test maps in \eqref{eq:sharpness-family} falls explicitly inside the unconstrained class $\mathcal A$. Calculating the perturbed energy validates the hypothesis that the energies remain bounded:
\begin{equation}\label{eq:sharpness-energy}
E^\eps(V^\eps)
=3\eps^{1/3}\int_\Omega|\nabla f|^2\dd
+\frac{1}{\eps} \int_\Omega (\eps^{1/2} f^3)^2\dd
\le 3\int_\Omega|\nabla f|^2\dd+\int_\Omega f^6\dd.
\end{equation}
For any exactly segregated competitor $W\in S$, the constraint $w_1 w_2 w_3 = 0$ requires that at almost every spatial coordinate $x$, at least one specific component $w_i(x)$ equals exactly zero. At this point, the coordinate-wise squared difference evaluates to
\[
|V^\eps(x)-W(x)|^2\ge (v_i^\eps(x) - 0)^2 = \eps^{1/3}f(x)^2.
\]
Integrating this point-wise minimum over the spatial domain we get
\[
\norm{V^\eps-W}_{L^2(\Omega)^3}
\ge\eps^{1/6}\norm{f}_{L^2(\Omega)}.
\]
The competitor constructed by 
\[
W^\eps=\eps^{1/6}(0,f,f)\in S,
\]
attains this equality, completing the proof of \eqref{eq:sharpness-distance}. If a bound holding an exponent $\beta>1/6$ were possible for this bounded family, it would imply:
\[
\norm{f}_{L^2(\Omega)}\le C\eps^{\beta-1/6}\to0,
\]
contradicting $f\not\equiv0$.
\end{proof}

\begin{remark}\label{rem:correction-rate}
Proposition~\ref{prop:sharpness} proves optimality of the exponent $1/6$ for general bounded-energy families. Its test family is not a family of almost-minimizers: for zero boundary data, $m_\eps=0$, whereas its penalty energy is $\int_\Omega f^6\dd>0$. The correction estimate does not measure error relative to a specified minimizer, and \eqref{eq:correction-energy} does not assert that the retraction decreases Dirichlet energy.
\end{remark}

\begin{corollary}\label{cor:little-o}
Under \eqref{eq:standing}, suppose
\[
V^\eps\in\mathcal A,\quad E^\eps(V^\eps)\le m_\eps+r_\eps,
\quad 0\le r_\eps\to0.
\]
Then
\begin{equation}\label{eq:little-o}
\|V^\eps-\calR(V^\eps)\|_{L^2(\Omega)^3}
=\sqrt3\,\dist_{L^2(\Omega)^3}(V^\eps,S)
=o(\eps^{1/6}).
\end{equation}
\end{corollary}

\begin{proof}
By \eqref{eq:whole-penalty} and the proof of Proposition~\ref{prop:correction},
\[
\eps^{-1/6}\|V^\eps-\calR(V^\eps)\|_{L^2}
\le\sqrt3\,|\Omega|^{1/3}
\left(\eps^{-1}\int_\Omega(v_1^\eps v_2^\eps v_3^\eps)^2\dd\right)^{1/6}
\to0.
\]
The equality of the two distances is \eqref{eq:exact-distance}.
\end{proof}
 Note that the proof yields only a little-$o$ statement,
$\|V^\eps-\calR(V^\eps)\|_{L^2}=o(\eps^{1/6})$: since the convergence
$\eps^{-1}\int_\Omega(v_1^\eps v_2^\eps v_3^\eps)^2\dd\to0$ comes without a rate,
it does not produce a bound of the form $C\eps^\beta$ for any fixed $\beta>1/6$.

\begin{corollary}\label{cor:smooth-correction}
Assume \eqref{eq:lip-data}. For every $V\in\mathcal A$ maintaining bounded perturbed energy $E^\eps(V)<\infty$, and taking an arbitrary tolerance $\tau>0$, there exists
\[
Z\in S\cap\Lip(\overline\Omega)^3\cap C^\infty(\Omega)^3
\]
such that
\begin{align}
\norm{Z-V}_{L^2(\Omega)^3}
&\le\sqrt3\,|\Omega|^{1/3}
\bigl(\eps E^\eps(V)\bigr)^{1/6}+\tau,\label{eq:smooth-correction-distance}\\
E^0(Z)&\le3E(V)+\tau.\label{eq:smooth-correction-energy}
\end{align}
\end{corollary}

\begin{proof}
Let $U \coloneqq \calR(V)$. By Proposition~\ref{prop:correction}, the correction $U$ belongs to the constrained class $S$ and satisfies the bounds
\begin{equation}\label{eq:proof-cor-bounds}
\norm{V - U}_{L^2(\Omega)^3} \le \sqrt{3}\,|\Omega|^{1/3}\bigl(\eps E^\eps(V)\bigr)^{1/6},
\qquad
E^0(U) \le 3E(V).
\end{equation}

Because $U \in S$ and the boundary data satisfy the Lipschitz condition \eqref{eq:lip-data}, we may apply Theorem~\ref{thm:main}. This theorem guarantees the existence of a sequence $\{U^n\}_{n=1}^\infty$ belonging to the regular class $S \cap \Lip(\overline\Omega)^3 \cap C^\infty(\Omega)^3$ such that
\[
U^n \longrightarrow U \quad \text{strongly in } H^1(\Omega)^3.
\]
Strong convergence in $H^1(\Omega)^3$ directly implies strong convergence in $L^2(\Omega)^3$ as well as convergence of the Dirichlet energies. Therefore,
\[
\lim_{n\to\infty} \norm{U^n - U}_{L^2(\Omega)^3} = 0, \qquad \lim_{n\to\infty} E^0(U^n) = E^0(U).
\]

Given the tolerance $\tau > 0$, we can choose a sufficiently large index $N$ such that the specific approximant $Z \coloneqq U^N$ satisfies both
\[
\norm{Z - U}_{L^2(\Omega)^3} \le \tau \qquad \text{and} \qquad E^0(Z) \le E^0(U) + \tau.
\]

To establish \eqref{eq:smooth-correction-distance}, we apply the triangle inequality in the $L^2$ norm and substitute the distance bound from \eqref{eq:proof-cor-bounds}:
\begin{align*}
\norm{Z - V}_{L^2(\Omega)^3} 
&\le \norm{Z - U}_{L^2(\Omega)^3} + \norm{U - V}_{L^2(\Omega)^3} \\
&\le \tau + \sqrt{3}\,|\Omega|^{1/3}\bigl(\eps E^\eps(V)\bigr)^{1/6}.
\end{align*}
The energy bound \eqref{eq:smooth-correction-energy} follows by using  \eqref{eq:proof-cor-bounds}. We get    
\[
E^0(Z) \le E^0(U) + \tau \le 3E(V) + \tau.
\]
This  completes the proof.
\end{proof}

\begin{corollary}[Exact approximations along minimizers]\label{cor:regular-minimizers}
Assume \eqref{eq:lip-data}. Let $U^{\eps_n}$ be penalized minimizers converging strongly in $H^1$ to $U\in S$ along the subsequence of Theorem~\ref{thm:minimizers}. There exist
\[
Z^n\in S\cap\Lip(\overline\Omega)^3\cap C^\infty(\Omega)^3
\]
such that
\begin{equation}\label{eq:regular-minimizer-convergence}
\norm{Z^n-U^{\eps_n}}_{H^1(\Omega)^3}\to0,
\quad Z^n\to U\quad\text{strongly in }H^1(\Omega)^3,
\end{equation}
and
\begin{equation}\label{eq:regular-minimizer-energy}
E^{\eps_n}(Z^n)=E^0(Z^n)\to m_0.
\end{equation}
\end{corollary}

\begin{proof}
By Theorem~\ref{thm:minimizers}, the sequence of penalized minimizers satisfies $U^{\eps_n} \to U$ strongly in $H^1(\Omega)^3$. By Lemma~\ref{lem:strong-continuity}, the retraction mapping $\calR$ is strongly continuous with respect to the $H^1$ topology. Therefore, applying the retraction to this sequence yields
\[
\calR(U^{\eps_n}) \longrightarrow \calR(U) \quad \text{strongly in } H^1(\Omega)^3.
\]
Since the target limit $U \in S$, it satisfies the partial segregation constraint, ensuring $\calR(U) = U$. 

For each index $n$, the retracted configuration $\calR(U^{\eps_n})$ belongs to $S$. Because the boundary data satisfy the Lipschitz condition \eqref{eq:lip-data}, we may apply Theorem~\ref{thm:main} to $\calR(U^{\eps_n})$. This guarantees the existence of a regular approximation. We select a specific approximant 
\[
Z^n \in S \cap \Lip(\overline\Omega)^3 \cap C^\infty(\Omega)^3
\]
satisfying the distance bound
\[
\norm{Z^n - \calR(U^{\eps_n})}_{H^1(\Omega)^3} < \frac{1}{n}.
\]
To establish strong convergence to $U$, we apply the triangle inequality in $H^1$ as follows
\[
\norm{Z^n - U}_{H^1(\Omega)^3} \le \norm{Z^n - \calR(U^{\eps_n})}_{H^1(\Omega)^3} + \norm{\calR(U^{\eps_n}) - U}_{H^1(\Omega)^3}.
\]
As $n \to \infty$, the first term is bounded by $1/n \to 0$, and the second term vanishes because $\calR(U^{\eps_n}) \to U$. Thus, $Z^n \to U$ strongly in $H^1(\Omega)^3$. 

Applying the triangle inequality again  gives
\[
\norm{Z^n - U^{\eps_n}}_{H^1(\Omega)^3} \le \norm{Z^n - U}_{H^1(\Omega)^3} + \norm{U - U^{\eps_n}}_{H^1(\Omega)^3}.
\]
Since both $Z^n \to U$ and $U^{\eps_n} \to U$ strongly in $H^1(\Omega)^3$, this difference vanishes, which verifies \eqref{eq:regular-minimizer-convergence}.

Finally, because $Z^n \in S$, its components satisfy $z_1^n z_2^n z_3^n = 0$ almost everywhere in $\Omega$. Consequently, the penalization integral evaluates to zero, yielding 
$$
E^{\eps_n}(Z^n) = E^0(Z^n) = E(Z^n).
$$ 
Since the Dirichlet energy is continuous with respect to strong $H^1$ convergence, $Z^n \to U$ implies
\[
\lim_{n\to\infty} E^{\eps_n}(Z^n) = \lim_{n\to\infty} E(Z^n) = E(U).
\]
By Theorem~\ref{thm:minimizers}, the limit $U$ is a global minimizer of $E^0$, so $E(U) = m_0$. This proves \eqref{eq:regular-minimizer-energy} and completes the proof.
\end{proof}

\section{Conclusion}

In summary, compatible $H^{1/2}$ boundary traces are sufficient to establish admissibility, the existence of minimizers, and the $\Gamma$-convergence of the penalized energies in arbitrary dimensions. By upgrading to Lipschitz boundary data, we obtain globally Lipschitz, interior-smooth approximations that strictly preserve both the exact boundary trace and the exact ternary segregation constraint. It is noteworthy that this constructive procedure requires no structural assumptions regarding finite interface decompositions, local junction profiles, or boundary phase collars.

The explicit geometric retraction serves a dual analytical purpose: it guarantees a continuous correction into the constrained class within the $H^1$ topology, and it provides a rigorous quantitative bound on the Dirichlet energy of the corrected map. Although the exact $L^2$-distance to the segregated space is smaller by a factor of $\sqrt{3}$, the pointwise nearest map may fail to have $H^1$ regularity, underscoring the necessity of the retraction mapping. For penalized almost-minimizers, the retraction error inherently vanishes in $H^1$ and decays at a rate of $o(\eps^{1/6})$ in $L^2$.

Finally, we emphasize that these smooth approximants are designed as variational competitors rather than exact solutions to the penalized Euler--Lagrange system. As a result, their respective positivity sets may naturally differ from the precise spatial supports of the original Sobolev configurations.

\bibliographystyle{acm}
\bibliography{refs}

@article{ContiTerraciniVerzini2005,
  author  = {Conti, Monica and Terracini, Susanna and Verzini, Gianmaria},
  title   = {Asymptotic estimates for the spatial segregation of competitive systems},
  journal = {Advances in Mathematics},
  volume  = {195},
  number  = {2},
  pages   = {524--560},
  year    = {2005},
  doi     = {10.1016/j.aim.2004.08.006}
}

@article{CaffarelliLin2007,
  author  = {Caffarelli, Luis A. and Lin, Fang-Hua},
  title   = {Singularly perturbed elliptic systems and multi-valued harmonic functions with free boundaries},
  journal = {Journal of the American Mathematical Society},
  volume  = {21},
  number  = {3},
  pages   = {847--862},
  year    = {2008},
  doi     = {10.1090/S0894-0347-08-00593-6}
}

@article{chang2004segregated,
  author  = {Chang, Shu-Ming and Lin, Chang-Shou and Lin, Tai-Chia and Lin, Wen-Wei},
  title   = {Segregated nodal domains of two-dimensional multispecies {B}ose--{E}instein condensates},
  journal = {Physica D: Nonlinear Phenomena},
  volume  = {196},
  number  = {3--4},
  pages   = {341--361},
  year    = {2004},
  doi     = {10.1016/j.physd.2004.06.002}
}

@article{tavares2003segregated,
  author  = {Tavares, Hugo and Terracini, Susanna},
  title   = {Regularity of the nodal set of segregated critical configurations under a weak reflection law},
  journal = {Calculus of Variations and Partial Differential Equations},
  volume  = {45},
  number  = {3--4},
  pages   = {273--317},
  year    = {2012},
  doi     = {10.1007/s00526-011-0458-z}
}

@article{ArakelyanBozorgnia2017,
  author  = {Arakelyan, Avetik and Bozorgnia, Farid},
  title   = {Uniqueness of limiting solution to a strongly competing system},
  journal = {Electronic Journal of Differential Equations},
  volume  = {2017},
  number  = {96},
  pages   = {1--8},
  year    = {2017},
  url     = {https://ejde.math.txstate.edu/Volumes/2017/96/abstr.html}
}

@article{CaffarelliPatriziQuitalo2017,
  author  = {Caffarelli, Luis A. and Patrizi, Stefania and Quitalo, Veronica},
  title   = {On a long range segregation model},
  journal = {Journal of the European Mathematical Society},
  volume  = {19},
  number  = {12},
  pages   = {3575--3628},
  year    = {2017},
  doi     = {10.4171/JEMS/747}
}

@article{Bozorgnia2016Acta,
  author  = {Bozorgnia, Farid},
  title   = {Uniqueness result for long range spatially segregation elliptic system},
  journal = {Acta Applicandae Mathematicae},
  volume  = {154},
  number  = {1},
  pages   = {1--14},
  year    = {2018},
  doi     = {10.1007/s10440-017-0129-y}
}

@article{bozorgnia2018singularly,
  author  = {Bozorgnia, Farid and Burger, Martin and Fotouhi, Morteza},
  title   = {On a class of singularly perturbed elliptic systems with asymptotic phase segregation},
  journal = {Discrete and Continuous Dynamical Systems},
  volume  = {42},
  number  = {7},
  pages   = {3539--3556},
  year    = {2022},
  doi     = {10.3934/dcds.2022023}
}

@article{soave2024partial1,
  title={On some singularly perturbed elliptic systems modeling partial segregation: uniform H$\backslash$" older estimates and basic properties of the limits},
  author={Soave, Nicola and Terracini, Susanna},
  journal={arXiv preprint arXiv:2409.11976},
  year={2024}
}

@article{soave2024partial2,
  title={On partially segregated harmonic maps: optimal regularity and structure of the free boundary},
  author={Soave, Nicola and Terracini, Susanna},
  journal={arXiv preprint arXiv:2410.23976},
  year={2024}
}

@book{braides2002gamma,
  author    = {Braides, Andrea},
  title     = {{$\Gamma$}-Convergence for Beginners},
  series    = {Oxford Lecture Series in Mathematics and Its Applications},
  volume    = {22},
  publisher = {Oxford University Press},
  address   = {Oxford},
  year      = {2002},
  isbn      = {9780198507840}
}

@article{ArakelyanBozorgnia2012,
  author  = {Arakelyan, Avetik and Bozorgnia, Farid},
  title   = {Numerical algorithms for a variational problem of the spatial segregation of reaction--diffusion systems},
  journal = {Applied Mathematics and Computation},
  volume  = {219},
  number  = {17},
  pages   = {8863--8875},
  year    = {2013},
  doi     = {10.1016/j.amc.2013.03.074}
}

@article{Bozorgnia2009,
  author  = {Bozorgnia, Farid},
  title   = {Numerical algorithm for spatial segregation of competitive systems},
  journal = {SIAM Journal on Scientific Computing},
  volume  = {31},
  number  = {5},
  pages   = {3946--3958},
  year    = {2009},
  doi     = {10.1137/080722588}
}

@article{arakelyan2016numerical,
  author  = {Arakelyan, Avetik and Barkhudaryan, Rafayel},
  title   = {A numerical approach for a general class of the spatial segregation of reaction--diffusion systems arising in population dynamics},
  journal = {Computers \& Mathematics with Applications},
  volume  = {72},
  number  = {11},
  pages   = {2823--2838},
  year    = {2016},
  doi     = {10.1016/j.camwa.2016.10.007}
}

@article{arakelyan2018convergence,
  author  = {Arakelyan, Avetik},
  title   = {Convergence of the finite difference scheme for a general class of the spatial segregation of reaction--diffusion systems},
  journal = {Computers \& Mathematics with Applications},
  volume  = {75},
  number  = {12},
  pages   = {4232--4240},
  year    = {2018},
  doi     = {10.1016/j.camwa.2018.03.025}
}

@misc{Vauchelet2025,
  author        = {Vauchelet, Nicolas},
  title         = {On a reaction--diffusion system modeling strong competition between two mosquito populations},
  year          = {2025},
  eprint        = {2506.15202},
  archivePrefix = {arXiv},
  primaryClass  = {math.AP},
  url           = {https://arxiv.org/abs/2506.15202}
}

@article{Xu2014,
  author  = {Xu, Benlong and Ni, Zhenzhang},
  title   = {Permanence of diffusive models for three competing species in heterogeneous environments},
  journal = {Abstract and Applied Analysis},
  volume  = {2014},
  pages   = {Article ID 376919, 10 pages},
  year    = {2014},
  doi     = {10.1155/2014/376919}
}

@article{LanzaraMontefusco2019,
  author  = {Lanzara, Flavia and Montefusco, Eugenio},
  title   = {On the limit configuration of four species strongly competing systems},
  journal = {Nonlinear Differential Equations and Applications NoDEA},
  volume  = {26},
  number  = {3},
  pages   = {Article 19, 17 pages},
  year    = {2019},
  doi     = {10.1007/s00030-019-0565-7}
}

@misc{bn,
  author        = {Bozorgnia, Farid and Arakelyan, Avetik and Kungurtsev, Vyacheslav and Valdman, Jan},
  title         = {Numerical Algorithms for Partially Segregated Elliptic Systems},
  year          = {2026},
  eprint        = {2603.05991},
  archivePrefix = {arXiv},
  primaryClass  = {math.NA},
  url           = {https://arxiv.org/abs/2603.05991}
}

@book{KinderlehrerStampacchia1980,
  author    = {Kinderlehrer, David and Stampacchia, Guido},
  title     = {An Introduction to Variational Inequalities and Their Applications},
  series    = {Pure and Applied Mathematics},
  volume    = {88},
  publisher = {Academic Press},
  address   = {New York},
  year      = {1980},
  isbn      = {978-0124073500}
}

@book{EvansGariepy2015,
  author    = {Evans, Lawrence C. and Gariepy, Ronald F.},
  title     = {Measure Theory and Fine Properties of Functions},
  edition   = {Revised},
  publisher = {CRC Press},
  address   = {Boca Raton, FL},
  year      = {2015},
  isbn      = {978-1482242386}
}

@article{bozorgnia2026interfaces,
  title={Interface distributions in partial segregation problems},
  author={Bozorgnia, Farid},
  journal={Applicable Analysis},
  pages={1--20},
  year={2026},
  publisher={Taylor \& Francis}
}

\end{document}